\documentclass[a4paper,10.9pt]{article}
\usepackage{bbding}
\usepackage{}
\usepackage{amsmath}
\usepackage{mathrsfs}
\usepackage{graphicx, subfigure}
\usepackage{wasysym}
\usepackage{color}
\usepackage{bm}
\usepackage{pifont}
\usepackage{txfonts}
\usepackage{graphicx}
\usepackage{amsfonts}
\usepackage{amssymb}
\usepackage{mathrsfs,psfrag,eepic,epsfig}

\makeatletter \@addtoreset{equation}{section}

\makeatletter
\newfont{\footsc}{cmcsc10 at 8truept}
\newfont{\footbf}{cmbx10 at 8truept}
\newfont{\footrm}{cmr10 at 10truept}
\title{\bf{Long-time behavior for the Cauchy problem of the 3-component Manakov system} \footnote{This work is supported by the National Key Research and Development Program of
China under Grant No. 2017YFB0202901 and the National Natural Science Foundation of China under Grant No.11871180.\protect\\
\hspace*{3ex} $^{\ddag}$Corresponding authors.\protect\\
\hspace*{3ex} E-mail address: xiubinwang@163.com (X.B. Wang) and bohan@hit.edu.cn (B. Han)}}

\author{Xiu-Bin Wang$^{\ddag}$ and Bo Han \\
\small Department of Mathematics, Harbin Institute of Technology, Harbin 150001,
People's Republic of China}

\date{}
\begin{document}
\maketitle

\noindent{\large \bf Abstract:}
In this work,
the Riemann-Hilbert problem for the 3-component Manakov system
is formulated on the basis of the corresponding $4\times 4$ matrix spectral problem.
Furthermore, by applying the nonlinear steepest descent techniques to an associated $4\times 4$ matrix valued Riemann-Hilbert problem,
we can find leading-order asymptotics for the Cauchy problem of the 3-component Manakov system.


\section{Introduction}

It is well known that distinct approaches \cite{vbma-1991}-\cite{tsff-cpaa} had been presented to analyze the nonlinear
problems, among which the inverse scattering transform (IST) \cite{cs-1967} had been successfully applied to
many important integrable systems. Particularly,
Zakharov and his collaborators \cite{ve-1984,ve-2} developed a Riemann-Hilbert formulation which is a modern version of the IST.
The Riemann-Hilbert problem \cite{xtc-2019}-\cite{wxb-2019}  is a nonlinear mapping between the set of smooth
potentials and the corresponding set of spectral data, which is a powerful tool for
the construction of soliton solutions of nonlinear evolution equations.
It is also possible
to investigate effective asymptotic results by performing the asymptotic analysis of the
corresponding Riemann-Hilbert problems.
To do asymptotics, one usually has to formulate a Riemann-Hilbert problem that has a unique solution that depends
on the scattering data.
At the beginning of the 1990's,
the nonlinear steepest descent technique (also called Deift-Zhou method) introduced by Deift
and Zhou \cite{PA-1993} presents a detailed rigorous proof to analyze the time asymptotic behaviors of the
nonlinear integrable partial differential equations (PDEs).
This technique was inspired by earlier work of Its \cite{Its-1973} and Manakov \cite{SV-1973};
for a detailed historical information, please see Ref.\cite{PAD-1993}, further developed by Deift, Venakides, and Zhou \cite{PS-1994}.
This technique has been successfully applied in determining asymptotic formulas for the initial value
problems of a large number of integrable systems associated with $2\times 2$ and $3\times 3$ Lax pairs
such as nonlinear Schr\"{o}dinger (NLS) equation \cite{XZJJ-4,LMP-1}, Korteweg-de Vries (KdV) equation \cite{KG-2009},
modified KdV equation \cite{PA-1993},
sine-Gordon (SG) equation \cite{PJ-1999}, Camassa-Holm (CH) equation \cite{bdma-2},
derivative NLS equation \cite{XJ-2013,tsf-2018}, Sasa-Satsuma (SS) equation \cite{gbl-jmp}, coupled NLS equation \cite{gxg-2018},
Degasperis-Procesi equation \cite{AB-2013},  modified NLS equation \cite{AV-1999},
three-component coupled modified KdV equation \cite{mwx-2019}, Fokas-Lenells equation \cite{XJ-2015}, Harry Dym equation \cite{LH-2015},
Kundu-Eckhaus equation \cite{wds-jde,gbl-jmp-2018}, and nonlocal NLS equation \cite{yr-jmp-1} etc.

%

In recent years, the investigation of multi-component NLS equations has been paid much attention,
since they can describe a variety of complex physical phenomena and
admit more abundant dynamics of localized wave solutions than ones in the scalar equations.
In this work, we therefore focus on the well-known 3-component Manakov equation \cite{gz-2018}-\cite{GDB-2015},
whose form yields
\begin{equation}\label{tc-NLS}
\left\{ \begin{aligned}
&iq_{1,t}+\frac{1}{2}q_{1,xx}+\left(|q_{1}|^2+|q_{2}|^2+|q_{3}|^2\right)q_{1}=0,\\
&iq_{2,t}+\frac{1}{2}q_{2,xx}+\left(|q_{1}|^2+|q_{2}|^2+|q_{3}|^2\right)q_{2}=0,\\
&iq_{3,t}+\frac{1}{2}q_{3,xx}+\left(|q_{1}|^2+|q_{2}|^2+|q_{3}|^2\right)q_{3}=0,~~(x,t)\in\mathbb{R}\times[0,+\infty),\\
&q_{1}(x,0)=q_{1,0}(x),~~q_{2}(x,0)=q_{2,0}(x),~~q_{3}(x,0)=q_{3,0}(x),
             \end{aligned} \right.
\end{equation}
where $q_{i}(x,t)(i=1,2,3)$ are complex-valued.
Besides, $q_{1}(x,0)$, $q_{2}(x,0)$,$q_{3}(x,0)$ lie in the Sobolev space
\begin{equation*}\label{NLS-4}
\mathcal {H}^{1,1}(\mathbb{R})=\left\{f(x)\in\mathscr{L}^{2}(\mathbb{R}):f'(x),xf(x)\in\mathscr{L}^{2}\right\}.
\end{equation*}
In addition, $q_{1,0}(x)$, $q_{2,0}(x)$ and $q_{3,0}(x)$ are assumed to be
generic so that $\det a(\lambda)$ defined in the following context is nonzero in $\mathbb{C}_{-}$.
The set of such generic functions is an open dense subset of $\mathcal {H}^{1,1}(\mathbb{R})$ \cite{rbp-1,GDB-2015},  which we denote by $\mathcal {G}$.
To our knowledge, much research work has been done for multi-component models \cite{ZLC-2017}-\cite{BF-2015}.
For example,
the initial-boundary value (IBV) problem for the 3-component Manakov system \eqref{tc-NLS} on a finite interval is investigated via the Fokas method \cite{ZLC-2009}.
Besides, in Ref.\cite{BF-2015},
the authors derived novel dark-bright soliton solutions for the 3-component Manakov system \eqref{tc-NLS}
and study the resulting soliton interactions.
Recently, 
Tian studied the IBV problems for the two-component NLS equation on
the interval and on the half-line \cite{tsf-jde,tsf-pa}.
However, the long-time asymptotics for the 3-component Manakov system \eqref{tc-NLS},
to the best knowledge of the authors, has never been reported up to now.

As we all know,
the Deift-Zhou method is a powerful technique to
analyze the long-time asymptotics for integrable nonlinear PDEs.
However, since Eq.\eqref{tc-NLS} admits
a $4\times 4$ matrix spectral problem, the RHP for Eq.\eqref{tc-NLS}
is rather complicated to derive.
The research in this work, to the best knowledge of the authors, has not been considered so far.
The main purpose of the present article is to analyze the long-time asymptotics for Eq.\eqref{tc-NLS}
by utilizing nonlinear steepest descent technique.

The structure of this article is given as follows.
In section 2, we convert the solution of the
Cauchy problem for Eq.\eqref{tc-NLS} to that of a matrix RHP.
In section 3, we transform the original
RHP to an appropriate expression and discuss the long-time asymptotic behavior of the solution of Eq.\eqref{tc-NLS}.
Finally,
some conclusions are presented in section 4.

\section{Riemann-Hilbert Problem}

The 3-component Manakov system is of course completely integrable,
it admits the result of the compatibility between the following linear differential equations \cite{gz-2018}
\begin{equation}\label{Lax-1}
\left\{ \begin{aligned}
&\psi_{x}=\left(i\lambda\sigma+\textbf{U}\right)\psi,\\
&\psi_{t}=\left(i\lambda^2\sigma+\textbf{V}\right)\psi,
                 \end{aligned} \right.
\end{equation}
with
\begin{align}\label{Lax-3}
&\sigma=\mbox{diag}(-1,1,1,1),~~
\textbf{U}(x,t)=\left(
                  \begin{array}{cccc}
                    0 & q_{1} & q_{2} & q_{2} \\
                    \bar{q}_{1} & 0 & 0 & 0 \\
                    \bar{q}_{2} & 0 & 0 & 0 \\
                    \bar{q}_{3} & 0 & 0 & 0 \\
                  \end{array}
                \right)
=i\left(
             \begin{array}{cc}
               0 & \textbf{q} \\
               \textbf{q}^{\dag} & \textbf{0} \\
             \end{array}
           \right),\notag\\
&\textbf{q}=\left(q_{1}, q_{2}, q_{3}\right),~~
\textbf{V}=\lambda \textbf{U}+i(\textbf{U}_{x}+\textbf{U}^2)\sigma/2,
\end{align}
where the superscript ``\dag'' means Hermitian conjugate of a matrix,
and $\Psi=\Psi(x,t,\lambda)$ is a column vector function of the spectral parameter $\lambda$.

In the following, by introducing a new matrix function
\begin{equation}\label{SA-6}
\psi(x,t;\lambda)=\mu(x,t;\lambda)e^{i\left(\lambda x+\lambda^2t\right)\sigma}.
\end{equation}
Then the spectral problem \eqref{Lax-1} gives
\begin{equation}\label{SA-7}
\left\{ \begin{aligned}
&\mu_{x}(x,t;\lambda)-i\lambda[\sigma,\mu(x,t;\lambda)]=\textbf{U}(x,t)\mu(x,t;\lambda),\\
&\mu_{t}(x,t;\lambda)-i\lambda^2[\sigma,\mu(x,t;\lambda)]=\textbf{V}(x,t;\lambda)\mu(x,t;\lambda).
                       \end{aligned} \right.
\end{equation}

We next present two eigenfunctions $\mu_{\pm}(x,t;\lambda)$ of $x$-part of \eqref{SA-7} by the following Volterra type integral equations
\begin{equation}\label{SA-9}
\mu_{\pm}=\mathcal {I}+\int_{\pm\infty}^{x} e^{i\lambda(x-\xi)\hat{\sigma}}[\textbf{U}(\xi,t)\mu_{\pm}(\xi,t;\lambda)]d\xi,
\end{equation}
where $\hat{\sigma}$ represents the operators which act on a $4\times 4$ matrix $\Omega$ by $\hat{\sigma}=[\sigma,\Omega]$.
Then $e^{\hat{\sigma}}=e^{\sigma}\Omega e^{\sigma}$.

It can therefore be shown that the functions $\mu_{\pm}$ are bounded and analytical for $\lambda\in \mathbb{C}$,
while can belong to
\begin{equation*}\label{SA-10}
\left\{ \begin{aligned}
&\mu_{+}:(\mathbb{C}_{+},\mathbb{C}_{+},\mathbb{C}_{+},\mathbb{C}_{-}),\\
&\mu_{-}:(\mathbb{C}_{-},\mathbb{C}_{-},\mathbb{C}_{-},\mathbb{C}_{+}),
                       \end{aligned} \right.
\end{equation*}
where $\mathbb{C}_{+}$ and $\mathbb{C}_{-}$ represent the upper and lower half complex $\lambda$-plane, respectively.
Here taking $i\lambda\sigma=\mbox{diag}(z_{1},z_{2},z_{3},z_{4})$, we then have
\begin{equation}\label{SA-11}
\left\{ \begin{aligned}
&\mathbb{C}_{+}=\left\{\lambda\in\mathbb{C}|\mbox{Re}z_{1}=\mbox{Re}z_{2}=\mbox{Re}z_{3}>\mbox{Re}z_{4}\right\},\\
&\mathbb{C}_{-}=\left\{\lambda\in\mathbb{C}|\mbox{Re}z_{1}=\mbox{Re}z_{2}=\mbox{Re}z_{3}<\mbox{Re}z_{4}\right\}.
                       \end{aligned} \right.
\end{equation}
The solutions of differential equation \eqref{SA-7} can be related by a matrix independent of $x$ and $t$. As a result
\begin{equation}\label{SA-12}
\mu_{-}(x,t;\lambda)=\mu_{+}(x,t;\lambda)e^{i\left(\lambda x+i\lambda^2t\right)\hat{\sigma}}s(\lambda).
\end{equation}
Evaluation at $t=0$ arrives at
\begin{equation*}\label{SA-13}
s(\lambda)=\lim_{x\rightarrow+\infty}e^{-i\lambda x\hat{\sigma}}\mu_{-}(x,0;\lambda),
\end{equation*}
i.e.,
\begin{equation*}\label{SA-13}
s(\lambda)=\mathcal {I}+\int_{-\infty}^{+\infty}e^{i\lambda x\hat{\sigma}}\left[\textbf{U}(x,0)\mu_{-}(x,0;\lambda)\right]dx.
\end{equation*}
The fact that $\textbf{U}$ is traceless together with Eq.\eqref{SA-9} indicates
\begin{equation*}\label{SA-14-1}
\det(\mu_{\pm}(x,t;\lambda))=1.
\end{equation*}
As a consequence, one can obtain
\begin{equation}\label{SA-14}
 \det(s(\lambda))=1.
\end{equation}
Additionally, 
we can know that
\begin{equation*}\label{SA-15}
\textbf{U}=-\textbf{U}^{\dag}.
\end{equation*}
Then from \eqref{Lax-1} we find
\begin{equation}\label{SA-17}
\psi_{x}^{A}(x,t;\lambda)=\left(i\lambda\sigma+\textbf{U}(x,t)\right)^{T}\psi^{A}(x,t;\lambda),
\end{equation}
with $\psi^{A}(x,t;\lambda)=(\psi^{-1}(x,t;\lambda))^{T}$, where the superscript `T' represents a matrix transpose.
As a result, we have
\begin{equation}\label{SA-181}
\psi^{\dag}(x,t;\bar{\lambda})=\psi^{-1}(x,t;\lambda).
\end{equation}
These relations mean that the eigenfunctions $\mu_{\pm}(x,t;\lambda)$ meet
\begin{equation}\label{SA-18}
\mu^{\dag}_{\pm}(x,t;\bar{\lambda})=\mu_{\pm}^{-1}(x,t;\lambda).
\end{equation}
To sum up, the matrix-valued function $s(\lambda)$ admits the following symmetry
\begin{equation}\label{SA-19}
s^{\dag}(\bar{\lambda})=s^{-1}(\lambda).
\end{equation}
In the next moment, without otherwise specified, by matrix blocking we rewrite the $4\times4$ matrix \textbf{A} as
\begin{equation*}\label{RHP-22}
\textbf{A}=\left(
    \begin{array}{cc}
      \textbf{A}_{11} & \textbf{A}_{12} \\
      \textbf{A}_{21} & \textbf{A}_{22} \\
    \end{array}
  \right),
\end{equation*}
where $\textbf{A}_{11}$ is scalar and $\textbf{A}_{22}$ is a $3\times 3$ matrix. It follows from \eqref{SA-12}-\eqref{SA-19} that
\begin{align}\label{RHH-1}
s_{11}^{\dag}(\bar{\lambda})=\det\left(s_{22}(\lambda)\right),~~
s_{21}^{\dag}(\bar{\lambda})=-s_{12}(\lambda)\mbox{adj}(s_{22}(\lambda)),
\end{align}
where $\mbox{adj}(\textbf{B})$ represents the adjoint matrix of matrix $\textbf{B}$.
Because of the above expression \eqref{RHH-1}, we can rewrite $s(\lambda)$ as
\begin{equation}\label{RHH-3}
s(\lambda)=\left(
             \begin{array}{cc}
               \det(a^{\dag}(\bar{\lambda}))  & b(\lambda) \\
               -\mbox{adj}(a^{\dag}(\bar{\lambda}))b^{\dag}(\bar{\lambda}) & a(\lambda) \\
             \end{array}
           \right),
\end{equation}
where $a(\lambda)$ is a $3\times3$ matrix-valued function and $b(\lambda)$ is a row vector-valued function.
It follows that $a(\lambda)$ and $b(\lambda)$ meet
\begin{equation}\label{RHH-5}
\left\{ \begin{aligned}
&a(\lambda)=\mathcal {I}+i\int_{-\infty}^{+\infty}q^{\dag}(\xi,0)\mu_{-,12}(\xi,0;\lambda)dx,\\
&b(\lambda)=i\int_{-\infty}^{+\infty}e^{2i\xi\lambda}q(\xi,0)\mu_{-,22}(\xi,0;\lambda)dx.
                       \end{aligned} \right.
\end{equation}
Obviously, it is not hard to see that $a(\lambda)$ is analytic in $\mathbb{C}_{+}$.

\noindent
\textbf{Theorem 2.1}.
Let $M(x,t;\lambda)$ be analytic for $\lambda \in \mathbb{C}\setminus \mathbb{R}$ and satisfy the following RHP
\begin{equation}\label{TH-1}
\left\{ \begin{aligned}
&M_{+}(x,t;\lambda)=M_{-}(x,t;\lambda)J(x,t;\lambda),~~\lambda\in\mathbb{R},\\
&M(x,t;\lambda)\rightarrow\mathcal {I},~~~~~~~~~~~~~~~~~~~~~~~~~~~~~~\lambda\rightarrow\infty,
                       \end{aligned} \right.
\end{equation}
where
\begin{equation}\label{TH-2}
J(x,t;\lambda)=\left(
                 \begin{array}{cc}
                   1+\gamma(\lambda)\gamma^{\dag}(\bar{\lambda})  & -e^{-it\Theta}\gamma^{\dag}(\bar{\lambda}) \\
                   -e^{it\Theta}\gamma(\lambda)  &  \mathcal {I}\\
                 \end{array}
               \right),
\end{equation}
and
\begin{equation}\label{TH-3}
M_{\pm}(x,t;\lambda)=\lim_{\epsilon\rightarrow 0^{+}}M(x,t;\lambda\pm i\epsilon),
~~\Theta(x,t;\lambda)=2\left(\frac{\lambda x}{t}+\lambda^2\right),~~\gamma(\lambda)=b(\lambda)a^{-1}(\lambda),
\end{equation}
with $\gamma(\lambda)\in\mathcal {H}^{1,1}(\mathbb{R})$ and $\sup_{\lambda\in\mathbb{R}}\gamma(\lambda)<\infty$.
Then the solution of this RHP exists and is unique.
Take
\begin{equation}\label{TH-4}
\textbf{q}(x,t)=2\lim_{\lambda\rightarrow\infty}(\lambda M(x,t;\lambda))_{12},
\end{equation}
which can tackle the Cauchy problem of the 3-component Manakov system \eqref{tc-NLS}.

\noindent
\textbf{Proof}.
The existence and uniqueness for the solution of the above RHP is a consequence of a ``vanishing lemma'' for
the associated RHP with the vanishing condition at infinity $M(\lambda) = O(1/\lambda) (\lambda\rightarrow\infty)$.
This result holds because the jump matrix $J(x,t;\lambda)$ is positive definite (also see Ref.\cite{MJ-2003}).

From the symmetries of jump matrix $J(x,t;\lambda)$, one gets that $M(x,t;\lambda)$
and $(M^{\dag})^{-1}(x,t;\lambda)$ suit the same RHP \eqref{TH-2}.
As the uniqueness for the solution of the RHP,
we can get
\begin{equation}\label{TH5}
M(x,t;\lambda)=\left(M^{\dag}\right)^{-1}(x,t;\lambda).
\end{equation}
Note that the asymptotic expansion of $M(x,t;\lambda)$
\begin{equation*}\label{RHP-16}
M(x,t;\lambda)=\mathcal {I}+\frac{M_{1}(x,t)}{\lambda}+\frac{M_{2}(x,t)}{\lambda^2}+\frac{M_{3}(x,t)}{\lambda^3}+\cdots,~~\lambda\rightarrow\infty,
\end{equation*}
which shows that $(M_{1})_{12}=-(M_{1})_{21}$.

We next derive that $\textbf{q}(x,t)$ determined by \eqref{TH-4} solves the equation \eqref{tc-NLS}
with the help of the standard arguments of the dressing method \cite{ve1,ve2}.
The chief idea of the dressing method is to derive two linear operators $\mathscr{L}$ and $\mathscr{N}$,
that is, (I) $\mathscr{L}\mathscr{N}$ and $\mathscr{N}\mathscr{L}$ satisfy the same jump condition as $M$.
(II) $\mathscr{L}\mathscr{N}$ and $\mathscr{N}\mathscr{L}$ are of $O(1/\lambda)$ as $\lambda\rightarrow\infty$.
Take
\begin{equation}\label{RHP-17}
\left\{ \begin{aligned}
&M_{x}=\mathscr{L}M+i\lambda[\sigma,M]+\textbf{U}M,\\
&M_{t}=\mathscr{N}M+i\lambda^2[\sigma,M]+\textbf{V}M,
                    \end{aligned} \right.
\end{equation}
where $\textbf{U}$ and $\textbf{V}$ are expressed by \eqref{Lax-3}.
A simple computation can lead to
\begin{equation*}\label{RHP-18}
\left(\mathscr{L}M\right)_{+}=\left(\mathscr{L}M\right)_{-}J,~~
\left(\mathscr{N}M\right)_{+}=\left(\mathscr{N}M\right)_{-}J.
\end{equation*}
If $\textbf{q}$ can be given by \eqref{TH-4}, the $\mathscr{L}M$ admits the homogeneous RHP
\begin{equation}\label{RHP-19}
\left\{ \begin{aligned}
&\left(\mathscr{L}M\right)_{+}=\left(\mathscr{L}M\right)_{-}J,\\
&\left(\mathscr{L}M\right)=O(1/\lambda),~~\lambda\rightarrow\infty,
                  \end{aligned} \right.
\end{equation}
which arrives at
\begin{equation}\label{TJJ1}
\mathscr{L}M=0.
\end{equation}
Furthermore, by comparing the coefficients of $O(1/\lambda)$ in the asymptotic expansion of \eqref{TJJ1}, one can have
\begin{equation*}\label{TJJ2}
M_{1}^{(O)}=-i\sigma\textbf{U}/2,~~\left[M_{1}^{(D)}\right]_{x}=i\sigma \textbf{U}^2/2,
\end{equation*}
where the superscripts ``(O)'' and ``(D)'' represent the off-diagonal and diagonal parts of block matrix, respectively.
As a consequence, $\mathscr{N}M$ meets the following homogeneous RHP
\begin{equation*}\label{RHP-20}
\left\{ \begin{aligned}
&\left(\mathscr{N}M\right)_{+}=\left(\mathscr{N}M\right)_{-}J,\\
&\left(\mathscr{N}M\right)=O(1/\lambda),~~\lambda\rightarrow\infty,
                  \end{aligned} \right.
\end{equation*}
which means that
\begin{equation}\label{RHP-21}
\mathscr{N}M=0.
\end{equation}
The compatibility condition of \eqref{RHP-19} and \eqref{RHP-21} gives the equation \eqref{tc-NLS}.
This indicates that $\textbf{q}(x,t)$ given by \eqref{TH-4} can solve the equation \eqref{tc-NLS}.

\section{Long-Time Asymptotics}
Inspired by earlier works of Deift and Zhou \cite{PA-1993},
the stationary point of $\Theta(\lambda)$ is given by $\lambda_{0}$, i.e.,
$\frac{d\Theta}{d\lambda}=0$, where $\lambda_{0}=-x/2t$,
as a result, $\Theta=2(\lambda^2-2\lambda_{0}\lambda)$.
Here we put our attention to physically interesting region $|\lambda_{0}|\leq C$,
in which $C$ is a constant.

\subsection{ Factorization of the Jump Matrix}

We now consider that the jump matrix admits two distinct factorizations
\begin{equation*}\label{TRH-1}
J=\left\{ \begin{aligned}
&\left(
    \begin{array}{cc}
       1 & e^{-it\Theta}\gamma(\lambda) \\
       0 & \mathcal {I} \\
    \end{array}
  \right)\left(
           \begin{array}{cc}
             1 & 0 \\
             -e^{-it\Theta}\gamma^{\dag}(\bar{\lambda}) & \mathcal {I} \\
           \end{array}
         \right),\\
&\left(
  \begin{array}{cc}
   1 & 0 \\
   -\frac{e^{-it\Theta}\gamma^{\dag}(\bar{\lambda})}{1+\gamma(\lambda)\gamma^{\dag}(\bar{\lambda})} &  \mathcal {I} \\
  \end{array}
\right)\left(
         \begin{array}{cc}
          1+\gamma(\lambda)\gamma^{\dag}(\bar{\lambda})  & 0 \\
           0 & \left(\mathcal {I}+\gamma^{\dag}(\bar{\lambda})\gamma(\lambda)\right)^{-1} \\
         \end{array}
       \right)\left(
                \begin{array}{cc}
                  1 & -\frac{e^{-it\Theta}\gamma(\lambda)}{1+\gamma(\lambda)\gamma^{\dag}(\bar{\lambda})} \\
                  0 & \mathcal {I} \\
                \end{array}
              \right).
                                \end{aligned} \right.
\end{equation*}
We next present a function $\delta(\lambda)$ as the solution of the matrix problem
\begin{equation}\label{TRH-2}
\left\{ \begin{aligned}
&\delta_{+}(\lambda)=\left(\mathcal {I}+\gamma^{\dag}\gamma\right)\delta_{-}(\lambda),~~|\lambda|<\lambda_{0},\\
&\delta(\lambda)\rightarrow \mathcal {I},~~~~~~~~~~~~~~~~~~~~~~~~~\lambda\rightarrow\infty.
                                \end{aligned} \right.
\end{equation}
As the jump matrix $\left(\mathcal {I}+\gamma^{\dag}\gamma\right)$ is positive definite,
the vanishing lemma arrives at the existence and uniqueness of the function $\delta(\lambda)$.
Additionally, we have
\begin{equation}\label{TRH-3}
\left\{ \begin{aligned}
&\det(\delta_{+}(\lambda))=\left(1+|\gamma|^{2}\right)\det(\delta_{-}(\lambda)),~~|\lambda|<\lambda_{0},\\
&\det(\delta(\lambda))\rightarrow 1,~~~~~~~~~~~~~~~~~~~~~~~~~~~~~~~~\lambda\rightarrow\infty.
                                \end{aligned} \right.
\end{equation}
From the positive definiteness of the jump matrix
$\mathcal {I}+\gamma^{\dag}\gamma$ and the vanishing lemma \cite{MJ-2003},
one can infer that $\delta$ exists and is unique.
A simple and direct calculation reveals that $\det\delta$ can be solved by the Plemelj formula \cite{MJ-2003}
\begin{equation}\label{DD-1}
\det(\delta(\lambda))=e^{\chi(\lambda)},
\end{equation}
where
\begin{equation*}\label{DD-2}
\chi(\lambda)=\frac{1}{2\pi i}\int_{-\infty}^{\lambda_{0}}
\left(\frac{\log\left(1+|\gamma(\xi)|^2\right)}{\xi-\lambda}\right)d\xi.
\end{equation*}
In practice, the above integral is singular as $\lambda\rightarrow\lambda_{0}$.
We therefore write the integral in the another form
\begin{align*}\label{DD-3}
\int_{-\infty}^{\lambda_{0}}
\left(\frac{\log\left(1+|\gamma(\xi)|^2\right)}{\xi-\lambda}\right)d\xi&=
\log\left(1+|\gamma(\lambda_{0})|^2\right)\log(\lambda-\lambda_{0}) \notag\\
&-\log\left(1+|\gamma(\lambda_{0})|^2\right)\log(\lambda-\lambda_{0}+1) \notag\\
&+\int_{-\infty}^{\lambda_{0}-1}
\left(\frac{\log\left(1+|\gamma(\xi)|^2\right)}{\xi-\lambda}\right)d\xi \notag\\
&+\int_{\lambda_{0}}^{\lambda_{0}-1}
\left(\frac{\log\left(1+|\gamma(\xi)|^2\right)-\log\left(1+|\gamma(\lambda_{0})|^2\right)}{\xi-\lambda}\right)d\xi,
\end{align*}
which can indicate that all the terms with the exception of the term $\log(\lambda-\lambda_{0})$ are analytic
for $\lambda$ in a neighborhood of $\lambda_{0}$. As a result, $\det\delta$ can be written as
\begin{equation*}\label{DD-4}
\det(\delta(\lambda))=(\lambda-\lambda_{0})^{i\nu}e^{\widetilde{\chi}(\lambda)},
\end{equation*}
where
\begin{equation}\label{TRH-5}
\left\{ \begin{aligned}
&\nu=-\frac{1}{2\pi}\log\left(1+\left|\gamma(\lambda_{0})\right|^{2}\right)<0,\\
&\widetilde{\chi}(\lambda_{0})=\frac{1}{2\pi i}\left[\int_{-\infty}^{\lambda_{0}-1}
\left(\frac{\log\left(1+|\gamma(\xi)|^2\right)}{\xi-\lambda}\right)d\xi\right.\\
&\left.~~~~~~~~~~~+\int_{\lambda_{0}}^{\lambda_{0}-1}
\left(\frac{\log\left(1+|\gamma(\xi)|^2\right)-\log\left(1+|\gamma(\lambda_{0})|^2\right)}{\xi-\lambda}\right)d\xi\right].
      \end{aligned} \right.
\end{equation}
In addition, for $|\lambda|<\lambda_{0}$, it follows from \eqref{TRH-2} that
\begin{equation*}\label{TRH-7}
\lim_{\epsilon\rightarrow0^{+}}\delta(\lambda- i\epsilon)=\left(\mathcal {I}+\gamma(\lambda)^{\dag}\gamma(\lambda)\right)^{-1}
\lim_{\epsilon\rightarrow0^{-}}\delta(\lambda+ i\epsilon).
\end{equation*}
If we set $f(\lambda)=\left(\delta^{\dag}(\bar{\lambda})\right)^{-1}$, we then get
\begin{equation*}\label{TRH-8}
f_{+}(\lambda)=\left(\mathcal {I}+\gamma^{\dag}(\lambda)\gamma(\lambda)\right)f_{-}(\lambda).
\end{equation*}
Thus, we know
\begin{equation}\label{TRH-91}
\left(\delta^{\dag}(\bar{\lambda})\right)^{-1}=\delta(\lambda).
\end{equation}
Putting \eqref{TRH-91} into \eqref{TRH-2} yields
\begin{align*}\label{TRH-92}
&\left|\delta_{-}(\lambda)\right|^{2}
=\left\{ \begin{aligned}
&2-\frac{|\gamma(\lambda)|^{2}}{1+|\gamma(\lambda)|^{2}},~~\lambda<\lambda_{0},\\
&2,~~~~~~~~~~~~~~~~~~~~~~~\lambda>\lambda_{0},
      \end{aligned} \right.\notag\\
&\left|\delta_{+}(\lambda)\right|^{2}
=\left\{ \begin{aligned}
&2+\frac{|\gamma(\lambda)|^{2}}{1+|\gamma(\lambda)|^{2}},~~\lambda<\lambda_{0},\\
&2,~~~~~~~~~~~~~~~~~~~~~~~\lambda>\lambda_{0},
      \end{aligned} \right.\notag\\
&|\det\delta_{-}(\lambda)|\leq 1,~~|\det\delta_{-}(\lambda)|\leq 1+|\gamma(\lambda)|^{2}<\infty,
\end{align*}
for fixed $\lambda\in\mathbb{R}$.
Similar to Ref.\cite{gxg-2018}, after a direct calculation, we can obtain
\begin{equation}\label{TRH-9}
|\delta(\lambda)|\leq\mbox{const}<\infty,~~|\det\left(\delta(\lambda)\right)|\leq\mbox{const}<\infty,
\end{equation}
for all $\lambda$, where we define $|\textbf{A}|=\sqrt{\left(\mbox{tr}\textbf{A}^{\dag}\textbf{A}\right)}$ for all \textbf{A}.
Then define
\begin{equation*}\label{TRH-10}
\Delta(\lambda)=\left(
                  \begin{array}{cc}
                    \det\left(\delta(\lambda)\right) & 0 \\
                    0 & \delta(\lambda)^{-1} \\
                  \end{array}
                \right).
\end{equation*}
and
\begin{equation*}\label{TRH-101}
\rho(\lambda)=\left\{ \begin{aligned}
&\frac{-\gamma(\lambda)}{1+\gamma(\lambda)\gamma^{\dag}(\bar{\lambda})},~~\lambda<\lambda_{0},\\
&\gamma(\lambda),~~~~~~~~~~~~~~~~~~\lambda>\lambda_{0},
      \end{aligned} \right.
\end{equation*}
Introduce
\begin{equation}\label{TRH-11}
M^{\Delta}(x,t;\lambda)=M(x,t;\lambda)\Delta^{-1}(x,t;\lambda),
\end{equation}
and reverse the orientation for $\lambda>\lambda_{0}$ as seen in Fig.1, then $M^{\Delta}$ admits the following RHP
\begin{equation}\label{TRH-12}
\left\{ \begin{aligned}
&M^{\Delta}_{+}(x,t;\lambda)=M_{-}^{\Delta}(x,t;\lambda)J^{\Delta}(x,t;\lambda),~~\lambda\in\mathbb{R},\\
&M^{\Delta}(x,t;\lambda)\rightarrow\mathcal {I},~~~~~~~~~~~~~~~~~~~~~~~~~~~~~~~\lambda\rightarrow\infty,
      \end{aligned} \right.
\end{equation}
where
\begin{align*}\label{TRH-13}
&J^{\Delta}=\left(
              \begin{array}{cc}
                1 & 0 \\
                \frac{e^{it\Theta}\delta_{-}^{-1}(\lambda)\rho^{\dag}(\bar{\lambda})}{\det\delta_{-}(\lambda)} & \mathcal {I} \\
              \end{array}
            \right)\left(
                     \begin{array}{cc}
                       1 & e^{-it\Theta}\det\delta_{+}(\lambda)\rho(\lambda)\delta_{+}(\lambda) \\
                       0 & \mathcal {I} \\
                     \end{array}
                   \right).
\end{align*}
In the next subsection, we deform the contour. It is convenient to analyze analytic
approximations of $\rho(\lambda)$.

$~~~~~~~~~~~~~~~~~~~$
{\rotatebox{0}{\includegraphics[width=8.0cm,height=6.5cm,angle=0]{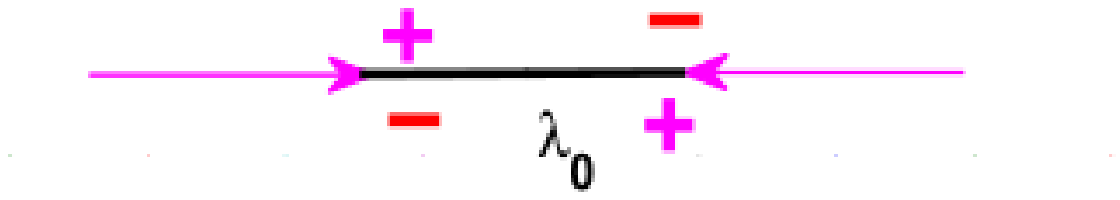}}}\\
\noindent {\small \textbf{Figure 1.}  The jump contour on $\mathbb{R}$ of \eqref{TRH-11}, and reverse the orientation for $\lambda>\lambda_{0}, \lambda_{0}=-x/2t$.\\

\subsection{Decomposition of $\rho(\lambda)$}

Let $\mathcal {L}$ denote the contour
\begin{equation}\label{DO-1}
\mathcal {L}:\left\{\lambda=\lambda_{0}+\alpha e^{\frac{3\pi i}{4}}:-\infty<\alpha<+\infty\right\}.
\end{equation}

\noindent
\textbf{Lemma 3.1} The vector-valued function $\rho(\lambda)$ has a decomposition
\begin{equation}\label{DO-2}
\rho(\lambda)=\mathcal {R}(\lambda)+\mathcal {H}_{1}(\lambda)+\mathcal {H}_{2}(\lambda),~~\lambda\in\mathbb{R},
\end{equation}
where $\mathcal {R}(\lambda)$ is piecewise rational and
$\mathcal {H}_{2}(\lambda)$ is analytically and continuously extended
to $\mathcal {L}$. In addition, $\mathcal {H}_{1}(\lambda)$ and $\mathcal {H}_{2}(\lambda)$ satisfy
\begin{equation}\label{AA-4}
\left\{ \begin{aligned}
&\left|e^{-it\Theta(\lambda)}\mathcal {H}_{1}(\lambda)\right|\lesssim \frac{1}{\left(1+|\lambda-\lambda_{0}|^{2}\right)t^{l}},~~\lambda\in\mathbb{R},\\
&\left|e^{-iT\theta(\lambda)}\mathcal {H}_{2}(\lambda)\right|\lesssim \frac{1}{\left(1+|\lambda-\lambda_{0}|^{2}\right)t^{l}},~~\lambda\in \mathcal {L},
      \end{aligned} \right.
\end{equation}
where positive integer $l$ is free. Taking the Schwartz conjugate
\begin{equation}\label{AA-5}
\rho^{\dag}(\bar{\lambda})=\mathcal {H}_{1}^{\dag}(\bar{\lambda})+\mathcal {H}_{2}^{\dag}(\bar{\lambda})+\mathcal {R}^{\dag}(\bar{\lambda}),
\end{equation}
leads to the same estimate for $e^{it\Theta(\lambda)}\mathcal {H}_{1}^{\dag}(\bar{\lambda})$,
$e^{it\Theta(\lambda)}\mathcal {H}_{2}^{\dag}(\bar{\lambda})$ and $e^{it\Theta(\lambda)}\mathcal {R}^{\dag}(\bar{\lambda})$
on the contour $\mathbb{R}\cup \bar{\mathcal {L}}$.

\noindent
\textbf{Proof}.
As $\lambda\geq\lambda_{0}$, $\rho(\lambda)=\gamma(\lambda)$.
With the help of the Taylor's expansion, we have
\begin{align*}\label{DO-4}
\left(\lambda-\lambda_{0}-i\right)^{m+5}\rho(\lambda)=
&\sum_{j=0}^{m}\rho_{j}\left(\lambda_{0}\right)\left(\lambda-\lambda_{0}\right)^{j}\notag\\
&+\frac{1}{m!}\int_{\lambda_{0}}^{\lambda}((\xi-\lambda_{0}-i)^{m+5}\rho(\xi))^{(m+1)}(\xi)(\lambda-\xi)^{m}d\xi,
\end{align*}
and define
\begin{equation*}\label{DO-5}
\mathcal {R}(\lambda)=\sum_{j=0}^{m}\rho_{j}(\lambda_{0})\left(\lambda-\lambda_{0}\right)^{j}/\left(\lambda-\lambda_{0}-i\right)^{m+5},~~
\mathcal {H}(\lambda)=\rho(\lambda)-\mathcal {R}(\lambda).
\end{equation*}
Then the following expression holds
\begin{equation*}\label{DO-6}
\frac{d^{j}\rho(\lambda)}{d\lambda^{j}}\big|_{\lambda=\lambda_{0}}=\frac{d^{j}\mathcal {R}(\lambda)}{d\lambda^{j}}\big|_{\lambda=\lambda_{0}},~~0\leq j\leq m.
\end{equation*}
For fixed $m\in\mathbb{Z}_{+}$, we next assume that m is of the form
\begin{equation*}\label{DO-7}
m=3p+1,~~p\in\mathbb{Z}_{+},
\end{equation*}
for convenience. As the map $\lambda\mapsto\Theta(\lambda)=2(\lambda^2-2\lambda_{0}\lambda)$ is one-to-one in $\lambda\geq\lambda_{0}$,
we define
\begin{equation}\label{DO-8}
f(\Theta)=\left\{ \begin{aligned}
&\frac{(\lambda-\lambda_{0}-i)^{p+2}}{(\lambda-\lambda_{0})^{p}}\mathcal {H}(\lambda)
=\frac{(\lambda-\lambda_{0})^{2p+2}}{(\lambda-\lambda_{0}-i)^{2p+4}},~~\Theta\geq-2\lambda_{0}^2,\\
&0,~~~~\Theta\leq-2\lambda_{0}^2,
      \end{aligned} \right.
\end{equation}
where
\begin{equation*}\label{DO-9}
g(\lambda)=\frac{1}{m!}\int_{0}^{1}\left((\cdot-\lambda_{0}-i)^{m+5}\rho(\cdot)\right)^{(m+1)}\left(\lambda_{0}
+\xi(\lambda-\lambda_{0})\right)(1-q_{1})^{m}d\xi,
\end{equation*}
which implies that
\begin{equation}\label{DO-10}
\left|\frac{d^{j}g(\lambda)}{d\lambda^{j}}\right|\lesssim 1,~~\lambda\geq\lambda_{0}.
\end{equation}
By making use of the Fourier transform, we have
\begin{equation*}\label{DO-11}
f(\Theta)=\frac{1}{\sqrt{2\pi}}\int_{-\infty}^{+\infty}e^{is\Theta}\hat{f}(s)ds,
\end{equation*}
where
\begin{equation*}\label{DO-12}
\hat{f}(s)=\frac{1}{\sqrt{2\pi}}\int_{-\infty}^{+\infty}e^{-is\Theta}f(\Theta)d\Theta,
\end{equation*}
Thus, as $0\leq j\leq p+1$,
\begin{equation*}\label{DO-13}
\int_{-2\lambda_{0}^2}\left|\frac{d^{j}f(\Theta)}{d\Theta^{j}}\right|^{2}d\Theta=
\int_{\lambda_{0}}^{+\infty}\left|\frac{(\lambda-\lambda_{0})^{2p+2-j}}{(\lambda-\lambda_{0}-i)^{2p+4}}\right|^{2}(\lambda-\lambda_{0})d\lambda\lesssim 1.
\end{equation*}
Resorting to the Plancherel theorem, we have
\begin{equation*}\label{DO-14}
\int_{-\infty}^{+\infty}\left(1+s^2\right)^{j}\left|\hat{f}(s)\right|^{2}ds\lesssim1,~~0\leq j\leq p+1.
\end{equation*}
We now present the decomposition of $\mathcal {H}(\lambda)$ as follows
\begin{align*}\label{DO-15}
&\mathcal {H}(\lambda)=\frac{(\lambda-\lambda_{0})^{p}}{\sqrt{2\pi}(\lambda-\lambda_{0}-i)^{p+2}}\int_{t}^{+\infty}e^{is\Theta}\hat{f}(s)ds
+\frac{(\lambda-\lambda_{0})^{p}}{\sqrt{2\pi}(\lambda-\lambda_{0}-i)^{p+2}}\int_{-\infty}^{t}e^{is\Theta}\hat{f}(s)ds \notag\\
&=\mathcal {H}_{1}(\lambda)+\mathcal {H}_{2}(\lambda).
\end{align*}
For $\lambda\geq\lambda_{0}$, we obtain
\begin{equation}\label{DO-16}
\left|e^{-it\Theta(\lambda)}\mathcal {H}_{1}(\lambda)\right|\lesssim\frac{|\lambda-\lambda_{0}-i|^{-2}}{\sqrt{2\pi}}\int_{t}^{+\infty}|\hat{f}(s)|ds
\lesssim\left|\lambda-\lambda_{0}-i\right|^{-2}t^{-p-\frac{1}{2}},
\end{equation}
where $0\leq r\leq p+1$. For $\lambda$ on the line $\{\lambda: \lambda_{0}+\epsilon e^{\frac{3\pi i}{4}}\}, \epsilon<0$, we have
\begin{equation}\label{DO-17}
\left|e^{-it\Theta(\lambda)}\mathcal {H}_{2}(\lambda)\right|\lesssim\left|\lambda-\lambda_{0}-i\right|^{-2}t^{-p/2},
\end{equation}
Here $m=3p+1$ is sufficiently large such that $p-1/2>p/2$
are greater than $l$. The proof finishes here.
Then the another case is also similar.

$~~~~~~~~~~~~~~~~~~~$
{\rotatebox{0}{\includegraphics[width=8.0cm,height=6.5cm,angle=0]{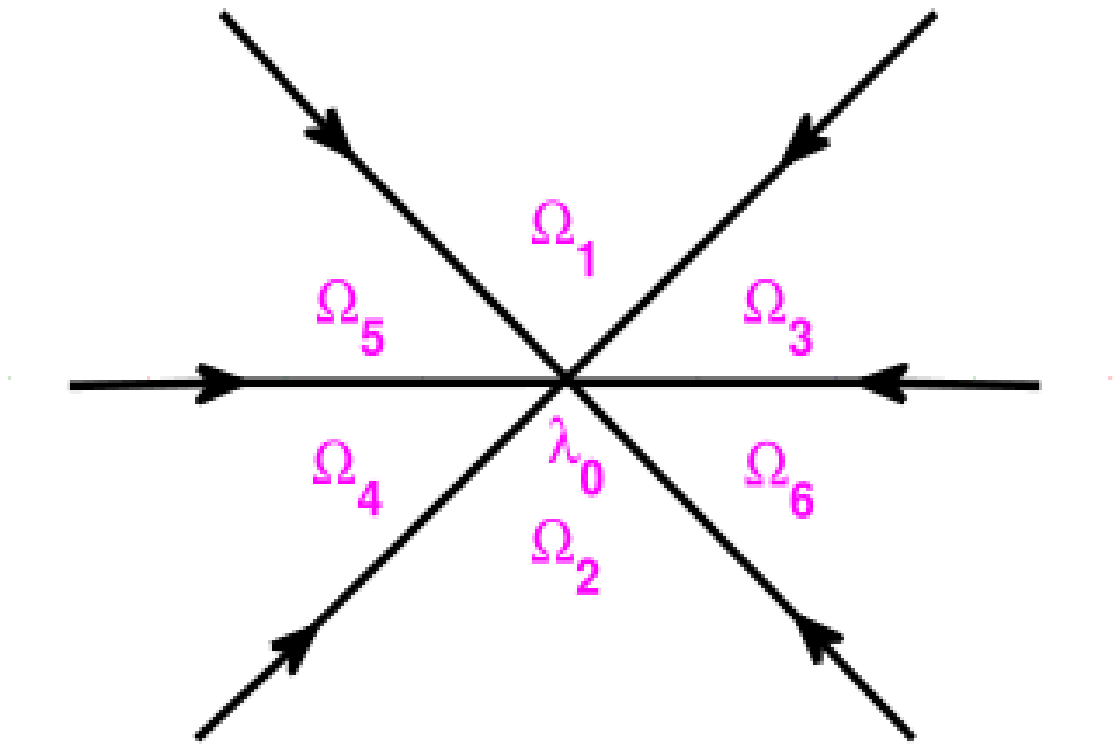}}}\\
\noindent {\small \textbf{Figure 2.} The orient jump contour
$\Sigma=\mathcal {L}\bigcup \bar{\mathcal {L}}\bigcup\mathbb{R}$ and domains $\{\Omega_{j}\}_{1}^{6}$,
and $\{\Omega_{j}\}_{1}^{6}$
denote the restriction of $\Sigma$ to the contour labeled by $j$. \\}

\subsection{First Contour Deformation}
In what follows, the original RHP turns into an equivalent RHP
formulated on an augmented contour $\Sigma$, where $\Sigma=\mathcal {L}\cup \bar{\mathcal {L}}\cup\mathbb{R}$ is shown in Fig.2.

Note that $J^{\Delta}(x,t;\lambda)$ can be rewritten as
\begin{equation*}\label{FC-1}
J^{\Delta}(x,t;\lambda)=(b_{-})^{-1}b_{+},
\end{equation*}
in which $b_{\pm}=\mathcal {I}\pm\omega_{\pm}$ with
\begin{equation*}\label{FC-2}
\left\{ \begin{aligned}
&b_{+}=b_{+}^{o}b_{+}^{a}=\left(\mathcal {I}+\omega_{+}^{o}\right)\left(\mathcal {I}+\omega_{+}^{a}\right)\\
&=\left(
   \begin{array}{cc}
     1 & e^{-it\Theta}\det \delta(\lambda)\mathcal {H}_{1}(\lambda)\delta(\lambda) \\
     0 & \mathcal {I}  \\
   \end{array}
 \right)\left(
          \begin{array}{cc}
            1 & e^{-it\Theta}\det \delta(\lambda)\left(\mathcal {H}_{2}(\lambda)+\mathcal {R}(\lambda)\right)\delta(\lambda) \\
            0 & \mathcal {I} \\
          \end{array}
        \right),\\
&b_{-}=b_{-}^{o}b_{-}^{a}=\left(\mathcal {I}-\omega_{-}^{o}\right)\left(\mathcal {I}-\omega_{-}^{a}\right)\\
&=\left(
   \begin{array}{cc}
      1 & 0 \\
     -\frac{e^{it\Theta}\delta^{-1}(\lambda)\mathcal {H}_{1}^{\dag}(\bar{\lambda})}{\det\delta(\lambda)} & \mathcal {I} \\
   \end{array}
 \right)\left(
          \begin{array}{cc}
            1 & 0 \\
            -\frac{e^{it\Theta}\delta^{-1}(\lambda)(\mathcal {H}_{2}^{\dag}(\bar{\lambda})+\mathcal {R}^{\dag}(\bar{\lambda}))}{\det\delta(\lambda)}) &  \mathcal {I} \\
          \end{array}
        \right).
              \end{aligned} \right.
\end{equation*}

\noindent
\textbf{Lemma 3.2} Define
\begin{equation}\label{CD-2}
M^{\sharp}(x,t;\lambda)=
\left\{ \begin{aligned}
&M^{\Delta}(x,t;\lambda),~~~~~~~~~~~~\lambda\in\Omega_{1}\cup\Omega_{2},\\
&M^{\Delta}(x,t;\lambda)(b^{a}_{-})^{-1},~~\lambda\in\Omega_{3}\cup\Omega_{4},\\
&M^{\Delta}(x,t;\lambda)(b^{a}_{+})^{-1},~~\lambda\in\Omega_{5}\cup\Omega_{6},
      \end{aligned} \right.
\end{equation}
As a consequence, the function $M^{\sharp}(\lambda)$ admits the RHP on the contour $\Sigma=\mathcal {L}\cup \bar{\mathcal {L}}\cup\mathbb{R}$ shown in Fig.2,
\begin{equation}\label{CD-3}
\left\{ \begin{aligned}
&M_{+}^{\sharp}(x,t;\lambda)=M_{-}^{\sharp}(x,t;\lambda)J^{\sharp}(x,t;\lambda),~~\lambda\in\Sigma,\\
&M^{\sharp}(x,t;\lambda)\rightarrow\mathcal {I},~~~~~~~~~~~~~~~\lambda\rightarrow\infty,
      \end{aligned} \right.
\end{equation}
where
\begin{equation}\label{CD-4}
J^{\sharp}=\left(b_{-}^{\sharp}\right)^{-1}b_{+}^{\sharp}=
\left\{ \begin{aligned}
&\left(b_{-}^{o}\right)^{-1}b_{+}^{o},~~\lambda\in\mathbb{R},\\
&\mathcal {I}^{-1}b_{+}^{a},~~~~~~~\lambda\in \mathbb{L},\\
&\left(b_{-}^{a}\right)^{-1}\mathcal {I},~~~\lambda\in \bar{\mathbb{L}}.
      \end{aligned} \right.
\end{equation}
We next construct the solution of the above RHP \eqref{CD-3} by using
the approach in Beals and Coifman \cite{rb-1984}. Assume that
\begin{equation*}\label{CD-5}
\omega^{\sharp}=\pm\left(b_{\pm}^{\sharp}-\mathcal {I}\right),~~\omega^{\sharp}=\omega^{\sharp}_{+}+\omega^{\sharp}_{-}.
\end{equation*}
Then denote the Cauchy operators $C_{\pm}$ for $\lambda\in\Sigma$ by
\begin{equation*}\label{CD-6}
C_{\pm}f(\lambda)=\frac{1}{2\pi i}\int_{\Sigma}\frac{f(\xi)}{\xi-\lambda_{\pm}}d\xi,
\end{equation*}
where $f\in\mathscr{L}^{2}(\Sigma)$. Define
\begin{equation}\label{CD-7}
C_{\omega^{\sharp}}f=C_{+}\left(f\omega_{-}^{\sharp}\right)+C_{-}\left(f\omega_{+}^{\sharp}\right).
\end{equation}

\noindent
\textbf{Theorem 3.3}. \cite{rb-1984} Assume $\mu^{\sharp}(x,t;\lambda)\in\mathscr{L}^{2}(\Sigma)+\mathscr{L}^{\infty}(\Sigma)$ satisfies
\begin{equation}\label{CD-8}
\mu^{\sharp}=\mathcal {I}+C_{\omega^{\sharp}}\mu^{\sharp}.
\end{equation}
Thus
\begin{equation}\label{CD-9}
M^{\sharp}(\lambda)=\mathcal {I}+\frac{1}{2\pi i}\int_{\Sigma}\frac{\mu^{\sharp}(\xi)\omega^{\sharp}(\xi)}{\xi-\lambda}d\xi,
\end{equation}
is the solution of the RHP \eqref{CD-3}.

\noindent
\textbf{Theorem 3.4}. The solutions $(q_{1}(x,t), q_{2}(x,t),q_{3}(x,t))$  of the Cauchy problem for the 3-component Manakov system \eqref{tc-NLS} can be expressed by
\begin{align}\label{CD-10}
\textbf{q}(x,t)&=(q_{1},q_{2},q_{3})
       =2\lim_{\lambda\rightarrow\infty}\left(\lambda M^{\sharp}(\lambda)\right)_{12}\notag\\
       &=\frac{i}{\pi}\left(\int_{\Sigma}\left(\mu^{\sharp}(x,t;\xi)\omega^{\sharp}(x,t;\xi)d\xi\right)d\xi\right)_{12}\notag\\
       &=\frac{i}{\pi}\left(\left(\int_{\Sigma}(1-C_{\omega^{\sharp}})^{-1}\mathcal {I}\right)(x,t;\xi)\omega^{\sharp}(x,t;\xi)d\xi\right)_{21}.
\end{align}
\textbf{Proof}.
It follows by \eqref{TRH-11}, \eqref{CD-3}, Theorem 2.1 and the following expressions
\begin{equation*}\label{CD-11}
\left\{ \begin{aligned}
&\left|e^{-it\Theta(\lambda)}\mathcal {H}_{2}(\lambda)\right|\lesssim|\lambda-\lambda_{0}-i|^{-2},~~\lambda\in\Omega_{5}\cup\Omega_{6},~~\lambda\rightarrow\infty,\\
&\left|e^{-it\Theta(\lambda)}\mathcal {R}(\lambda)\right|\lesssim|\lambda-\lambda_{0}-i|^{-5},~~~~\lambda\in\Omega_{5}\cup\Omega_{6},~~~\lambda\rightarrow\infty,.
      \end{aligned} \right.
\end{equation*}

\subsection{Second Contour Deformation}
We next reduce the RHP $M^{\sharp}$ on $\Sigma$ to the RHP $M'$ on $\Sigma'$,
where $\Sigma'=\Sigma \setminus \mathbb{R}=\mathcal {L}\cup \bar{\mathcal {L}}$ is orientated as in Fig.2.
In addition, we estimate the
error between the RHP on $\Sigma$ and $\Sigma'$. Then it is proved that the solution of
the equation \eqref{tc-NLS} can be expressed in terms of $M'$ by adding an error term.

\noindent
Let $\omega^{e}=\omega^{a}+\omega^{b}$, where\\
(I) $\omega^{a}=\omega^{\sharp}\upharpoonright \mathbb{R}$ and
composed of terms of type $\mathcal {H}_{1}(\lambda)$ and $\mathcal {H}_{1}^{\dag}(\bar{\lambda})$;\\
(II) $\omega^{b}$ is supported on $\mathcal {L}\cup \bar{\mathcal {L}}$  and composed of the contribution to $\omega^{\sharp}$ from terms of
type $\mathcal {H}_{2}(\lambda)$ and $\mathcal {H}_{2}^{\dag}(\bar{\lambda})$.

\noindent
Define $\omega'=\omega^{\sharp}-\omega^{e}$, then $\omega'=0$ on $\mathbb{R}$.
Thus, $\omega'$ is supported on $\Sigma'$ with contribution to $\omega^{\sharp}$ from $\mathcal {R}(\lambda)$ and $\mathcal {R}^{\dag}(\bar{\lambda})$.

\noindent
\textbf{Lemma 3.5}.
For sufficiently small $\epsilon$, as $t\rightarrow\infty$
\begin{equation}\label{CT-2}
\left\{ \begin{aligned}
&\left\|\omega^{a}\right\|_{\mathscr{L}^{\infty}(\mathbb{R})\bigcap\mathscr{L}^{1}(\mathbb{R})\bigcap\mathscr{L}^{2}(\mathbb{R})}\lesssim t^{-l},\\
&\left\|\omega^{b}\right\|_{{\mathscr{L}}^{\infty}(\mathcal {L}\bigcup \bar{\mathcal {L}})\bigcap\mathscr{L}^{1}(\mathcal {L}\bigcup \bar{\mathcal {L}})\bigcap\mathscr{L}^{2}(\mathcal {L}\bigcup \bar{\mathcal {L}})}\lesssim t^{-l},\\
&\|\omega'\|_{\mathscr{L}^{2}(\Sigma)}\lesssim t^{-1/4},~~\|\omega'\|_{\mathscr{L}^{1}(\Sigma)}\lesssim t^{-1/2}.
      \end{aligned} \right.
\end{equation}

\noindent
\textbf{Proof}.
The proof of first two expressions in \eqref{CT-2} follows from Lemma 2.2.
From the definition of $\mathcal {R}(\lambda)$, we have
\begin{equation*}\label{CT-3}
|\mathcal {R}(\lambda)|\lesssim\left(1+|\lambda-\lambda_{0}|^5\right)^{-1},
\end{equation*}
on the contour $\{\lambda=\lambda_{0}+\lambda_{0}\alpha e^{3\pi i/4}:-\infty<\alpha<+\infty\}$.
By means of inequality \eqref{TRH-9},
we obtain
\begin{equation}\label{CT-4}
\left|e^{-it\Theta}\det(\delta(\lambda))\delta(\lambda)\mathcal {R}(\lambda)\right|\lesssim e^{-2\alpha^2t}\left(1+|\lambda-\lambda_{0}|^5\right)^{-1},
\end{equation}
After a direct calculation, we can obtain the last expression in \eqref{CT-2}.

\noindent
\textbf{Lemma 3.6.}
In the case $0<\lambda_{0}\leq C$, as $t\rightarrow\infty$,
the inverse $(1-C_{\omega'})^{-1}:\mathscr{L}^{2}(\Sigma)\rightarrow\mathscr{L}^{2}(\Sigma)$
exists, and has uniform boundedness
\begin{equation*}\label{CT-5}
\left\|\left(1-C_{\omega'}\right)^{-1}\right\|_{\mathscr{L}^{2}(\Sigma)}\lesssim 1.
\end{equation*}
Besides
\begin{equation*}\label{CT-5}
\left\|\left(1-C_{\omega^{\sharp}}\right)^{-1}\right\|_{\mathscr{L}^{2}(\Sigma)}\lesssim 1.
\end{equation*}
\textbf{Proof.} See \cite{PA-1993} and references therein.

\noindent
\textbf{Lemma 3.7}. The integral equation has estimate as $t\rightarrow\infty$
\begin{equation}\label{CT-6}
\int_{\Sigma}\left(\left(1-C_{\omega^{\sharp}}\right)^{-1}\mathcal {I}\right)(\xi)\omega^{\sharp}(\xi)d\xi
=\int_{\Sigma}\left(\left(1-C_{\omega'}\right)^{-1}\mathcal {I}\right)(\xi)\omega'(\xi)d\xi+O\left(t^{-l}\right).
\end{equation}
\textbf{Proof.} Via a simple calculation, we derive
\begin{align*}\label{CT-7}
\left(\left(1-C_{\omega^{\sharp}}\right)^{-1}\mathcal {I}\right)\omega^{\sharp}
&=\left(\left(1-C_{\omega'}\right)^{-1}\mathcal {I}\right)\omega'
+\omega^{e}+\left(\left(1-C_{\omega'}\right)^{-1}\left(C_{\omega^{e}}\mathcal {I}\right)\right)\omega^{\sharp}\notag\\
&+\left(\left(1-C_{\omega'}\right)^{-1}\left(C_{\omega^{'}}\mathcal {I}\right)\right)\omega^{e}\notag\\
&+\left(\left(1-C_{\omega'}\right)^{-1}C_{\omega^{e}}\left(1-C_{\omega^{\sharp}}\right)^{-1}\right)\left(C_{\omega^{\sharp}}\mathcal {I}\right)\omega^{\sharp}.
\end{align*}
It follows from Lemma 3.5 that
\begin{equation*}\label{LEM-1}
\left\{ \begin{aligned}
&\left\|\omega^{e}\right\|_{\mathscr{L}^{1}(\Sigma)}\leq\|\omega^{a}\|_{\mathscr{L}^{1}(\mathbb{R})}+\|\omega^{b}\|_{\mathscr{L}^{1}(\mathcal {L}\cup\bar{\mathcal {L}})}\lesssim t^{-l},\\
&\|\left(\left(1-C_{\omega'})^{-1}(C_{\omega^{e}}\mathcal {I}\right)\right)\omega^{\sharp}\|_{\mathscr{L}^{1}(\Sigma)}\leq
\|\left(1-C_{\omega'}\right)^{-1}\|_{\mathscr{L}^{2}(\Sigma)}\|C_{\omega^{e}}\mathcal {I}\|_{\mathscr{L}^{2}(\Sigma)}\|\omega^{\sharp}\|_{\mathscr{L}^{2}(\Sigma)}\\
&~~~~~~~~~~\leq\|\omega^{e}\|_{\mathscr{L}^{2}(\Sigma)}\|\omega^{\sharp}\|_{\mathscr{L}^{2}(\Sigma)}\lesssim t^{-l-1/4},\\
&\|\left(\left(1-C_{\omega'}\right)^{-1}(C_{\omega^{'}}\mathcal {I})\right)\omega^{e}\|_{\mathscr{L}^{1}(\Sigma)}\leq
\left\|\left(1-C_{\omega'}\right)^{-1}\right\|_{\mathscr{L}^{2}(\Sigma)}\|C_{\omega^{'}}\mathcal {I}\|_{\mathscr{L}^{2}(\Sigma)}\|\omega^{e}\|_{\mathscr{L}^{2}(\Sigma)}\\
&~~~~~~~~~~\leq\|\omega^{'}\left\|_{\mathscr{L}^{2}(\Sigma)}\right\|\omega^{e}\|_{\mathscr{L}^{2}(\Sigma)}\lesssim t^{-l-1/4},\\
&\|((1-C_{\omega'})^{-1}C_{\omega^{e}}(1-C_{\omega^{\sharp}})^{-1})(C_{\omega^{\sharp}}\mathcal {I})\omega^{\sharp}\|_{\mathscr{L}^{1}(\Sigma)}\\
&~~~~~~~~~~=\|(1-C_{\omega'})^{-1}\|_{\mathscr{L}^{2}(\Sigma)}\|C_{\omega^{e}}\|_{\mathscr{L}^{2}(\Sigma)}\|(1-C_{\omega^{\sharp}})^{-1}\|_{\mathscr{L}^{2}(\Sigma)}
\|C_{\omega^{\sharp}}\mathcal {I}\|_{\mathscr{L}^{2}(\Sigma)}\|\omega^{\sharp}\|_{\mathscr{L}^{2}(\Sigma)}\\
&~~~~~~~~~~\lesssim\|\omega^{e}\|_{\mathscr{L}^{\infty}(\Sigma)}\|\omega^{\sharp}\|^{2}_{\mathscr{L}^{2}(\Sigma)}\lesssim t^{-l-1/2}.
      \end{aligned} \right.
\end{equation*}
This finishes the proof of the theorem.

\noindent
\textbf{Lemma 3.8} As $t\rightarrow\infty$
\begin{equation*}\label{LE10-1}
\textbf{q}(x,t)=\frac{i}{\pi}\left(\int_{\Sigma'}\left(\left(1-C_{\omega'}\right)^{-1}\mathcal {I}\right)(\xi)\omega'(\xi)d\xi\right)_{12}+O\left(t^{-l}\right).
\end{equation*}
\textbf{Proof} A straightforward consequence of Theorem 3.4 and Lemma 3.8.

\noindent
\textbf{Corollary 3.9} As $t\rightarrow\infty$,
\begin{equation}\label{CO11-1}
\textbf{q}(x,t)=2\lim_{\lambda\rightarrow\infty}\left(\lambda M'(x,t;\lambda)\right)_{12}+O\left(t^{-l}\right),
\end{equation}
where $M'(x,t;\lambda)$ satisfies the RHP
\begin{equation*}\label{CO11-2}
\left\{ \begin{aligned}
&M'_{+}(x,t;\lambda)=M'_{-}(x,t;\lambda)J'(x,t;\lambda),~~\lambda\in\Sigma',\\
&M'(x,t;\lambda)\rightarrow\mathcal {I},~~~~~~~~~~~~~~~~~~~~~~~~~~~\lambda\rightarrow\infty,
      \end{aligned} \right.
\end{equation*}
where
\begin{align*}\label{CO11-3}
&\omega'=\omega'_{+}+\omega'_{-},~~b_{\pm}'=\mathcal {I}+\omega_{\pm}',~~J'=(b'_{-})^{-1}b'_{+},\notag\\
&b'_{+}=\left(
         \begin{array}{cc}
           1 & e^{-it\Theta}\det(\delta(\lambda))\mathcal {R}(\lambda)\delta(\lambda) \\
           0 & \mathcal {I} \\
         \end{array}
       \right),~~b'_{-}=\mathcal {I},~~\mbox{on}~~\mathcal {L},\notag\\
&b'_{+}=\mathcal {I},~~b'_{-}=\left(
                               \begin{array}{cc}
                                 1 & 0 \\
                                 -\frac{e^{it\Theta}\delta^{-1}(\lambda)\mathcal {R}^{\dag}(\lambda)}{\det(\delta(\lambda))} & \mathcal {I}\\
                               \end{array}
                             \right),~~\mbox{on}~~\bar{\mathcal {L}},
\end{align*}

\noindent
\textbf{Proof}
Set $\mu'=(1-C_{\omega'})^{-1}\mathcal {I}$ and
\begin{equation*}\label{CO11-5}
M'(x,t;\lambda)=\mathcal {I}+\frac{1}{2\pi i}\int_{\Sigma'}\frac{\mu'(x,t;\xi)\omega'(x,t;\xi)}{\xi-\lambda}d\xi.
\end{equation*}
Similar to Theorem 3.4, we can construct this corollary 3.9 in terms of \eqref{LE10-1}.

$~~~~~~~~~~~~~~~~~~~$
{\rotatebox{0}{\includegraphics[width=8.0cm,height=6.5cm,angle=0]{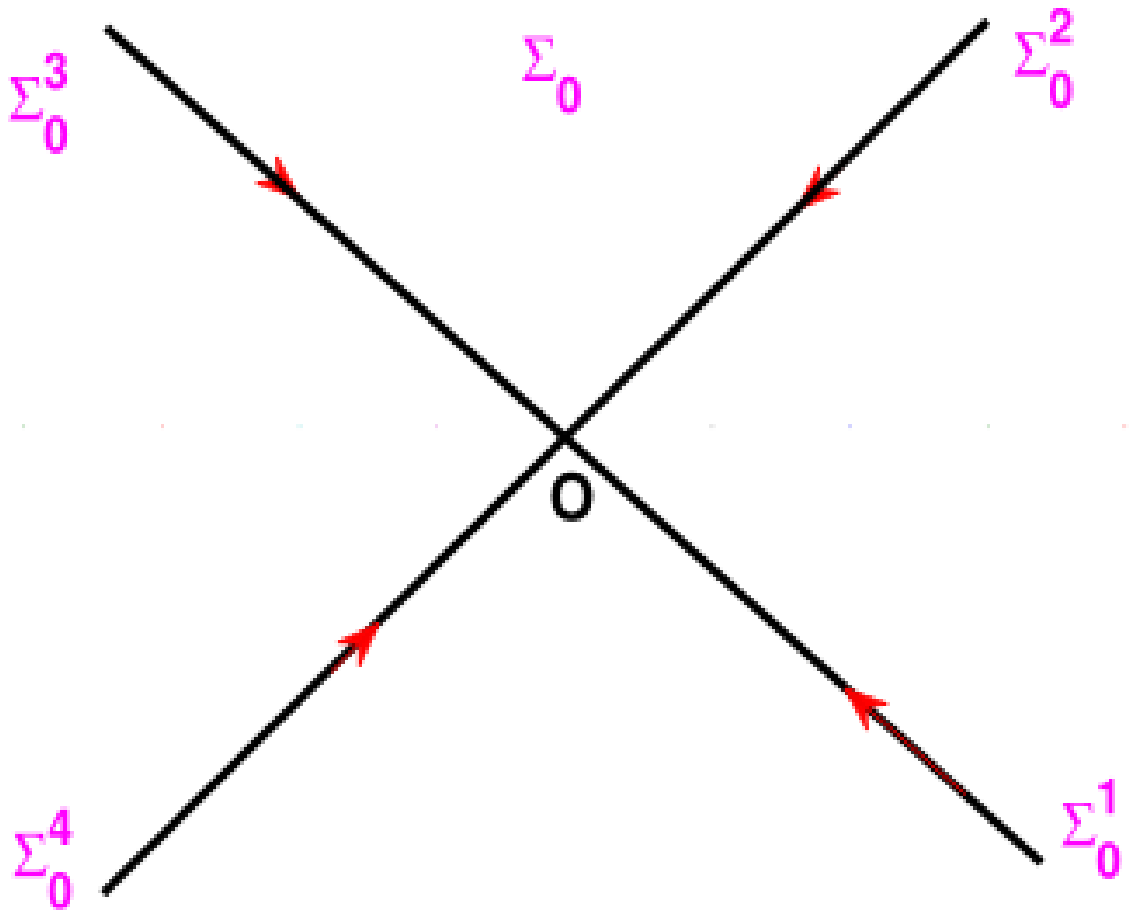}}}\\
\noindent {\small \textbf{Figure 3.} The oriented contour $\Sigma_{0}$ ($\Sigma_{0}$ denote the contour $\{\lambda=\alpha e^{\pm\frac{i\pi}{4}}:\alpha\in\mathbb{R}\}$). \\}

\subsection{Reduction of the Riemann-Hilbert Problems}
In this subsection, we localize the jump matrix of the RHP to the neighborhood
of the stationary phase point $\lambda_{0}$.
Under suitable scaling of the spectral parameter,
the RHP is reduced to a RHP with constant jump matrix which can be solved explicitly.

Let $\Sigma_{0}$ denote the contour $\{\lambda=\alpha e^{\pm\frac{i\pi}{4}}:\alpha\in\mathbb{R}\}$
oriented as in Fig.3. Define the scaling operator
\begin{align}\label{FR-1}
N:&\mathscr{L}^{2}(\Sigma')\rightarrow\mathscr{L}^{2}(\Sigma_{0}),\notag\\
&f(\lambda)\mapsto(Nf)(\lambda)=f\left(\lambda_{0}+\frac{\lambda}{2\sqrt{t}}\right),
\end{align}
and set $\hat{\omega}=N\omega'$. A direct change-of-variable argument means that
\begin{equation*}\label{FR-2}
C_{\omega'}=N^{-1}C_{\hat{\omega}}N,
\end{equation*}
where the operator $C_{\omega'}$ is a bounded map from $\mathscr{L}^{2}(\Sigma_{0})$ into $\mathscr{L}^{2}(\Sigma_{0})$.

One can infer that
\begin{equation*}\label{FR-3}
\hat{\omega}=\hat{\omega}_{+}=\left(
                                \begin{array}{cc}
                                  0 & (Ns_{1})(\lambda) \\
                                  0 & 0 \\
                                \end{array}
                              \right),
\end{equation*}
on $\hat{\mathcal {L}}=\{\lambda=\alpha e^{\frac{3\pi i}{4}}:-\infty<\alpha<+\infty\}$, and
\begin{equation*}\label{FR-4}
\hat{\omega}=\hat{\omega}_{-}=\left(
                                \begin{array}{cc}
                                  0 & 0\\
                                  (Ns_{2})(\lambda) & 0 \\
                                \end{array}
                              \right),
\end{equation*}
on $\hat{\bar{\mathcal {L}}}$, where
\begin{equation*}\label{FR-5}
s_{1}(\lambda)=e^{-it\Theta(\lambda)}(\det\delta)(\lambda\mathcal {R})(\lambda)\delta(\lambda),
~~s_{2}(\lambda)=\frac{e^{it\Theta(\lambda)}\delta^{-1}(\lambda)\mathcal {R}^{\dag}(\bar{\lambda})}{\det(\delta(\lambda))}.
\end{equation*}

\noindent
\textbf{Lemma 3.10} As $t\rightarrow\infty$, and $\lambda\in \hat{\mathcal {L}}$, for an arbitrary positive integer $l$, then
\begin{equation}\label{FR-6}
\left|(N\widetilde{\delta})(\lambda)\right|\lesssim t^{-l},
\end{equation}
where
\begin{equation}\label{FR-7}
\widetilde{\delta}(\lambda)=e^{-it\Theta(\lambda)}\mathcal {R}(\lambda)\left[\delta(\lambda)-(\det\delta)(\lambda)\mathcal {I}\right].
\end{equation}

\noindent
\textbf{Proof} It follows from \eqref{TRH-2} and \eqref{TRH-3}  that $\widetilde{\delta}$ satisfies the following RHP
\begin{equation}\label{FR-8}
\left\{ \begin{aligned}
&\widetilde{\delta}_{+}(\lambda)=\delta_{-}(\lambda)\left((1+|\gamma(\lambda)|^{2}\right)+e^{-it\Theta(\lambda)}f(\lambda),
~~\lambda\in(-\infty,\lambda_{0}),\\
&\widetilde{\delta}(\lambda)\rightarrow 0,~~~~~~~~~~~~~~~~~~~\lambda\rightarrow\infty,
      \end{aligned} \right.
\end{equation}
where $f(\lambda)=[\mathcal {R}(\gamma^{\dag}\gamma-|\gamma|^{2})\delta_{-}](\lambda)$.
The solution for the above RHP can be expressed by
\begin{align}\label{FR-10}
&\widetilde{\delta}(\lambda)=X(\lambda)\int_{-\infty}^{\lambda_{0}}\frac{e^{-it\Theta(\lambda)}f(\xi)}{X_{+}(\xi)(\xi-\lambda)}d\xi,\notag\\
&X(\lambda)=\exp\left(\frac{1}{2\pi i}\int_{-\infty}^{\lambda_{0}}\frac{\log\left(1+|\gamma(\lambda)|^2\right)}{\xi-\lambda}d\xi\right).
\end{align}
Observe that
\begin{equation}\label{FR-11}
\mathcal {R}\gamma^{\dag}\gamma-|\gamma|^{2}\mathcal {R}=(\mathcal {R}-\rho)\gamma^{\dag}\gamma-|\gamma|^{2}(\mathcal {R}-\rho).
\end{equation}
Similar to Lemma 3.2, $f(\lambda)$ can be decomposed into two parts: $f_{1}(\lambda)$ and $f_{2}(\lambda)$,
where $f_{2}(\lambda)$ admits an analytic continuation to $\mathcal {L}_{t}$ satisfying
\begin{equation}\label{FR-12}
\left\{ \begin{aligned}
&\left|e^{-it\Theta(\lambda)}f_{1}(\lambda)\right|\lesssim\frac{1}{\left(1+|\lambda-\lambda_{0}+1/t|^{2}\right)t^{l}},~~\lambda\in\mathbb{R},\\
&\left|e^{-it\Theta(\lambda)}f_{2}(\lambda)\right|\lesssim\frac{1}{\left(1+|\lambda-\lambda_{0}+1/t|^{2}\right)t^{l}},~~\lambda\in \mathcal {L}_{t},
      \end{aligned} \right.
\end{equation}
where (see Fig.4)
\begin{equation}\label{FR-13}
\mathcal {L}_{t}:\left\{\lambda=\lambda_{0}-1/t+\alpha e^{\frac{3\pi i}{4}}:0\leq\alpha<+\infty\right\}.
\end{equation}
As $\lambda\in\hat{\mathcal {L}}$, we have
\begin{align}\label{FR-14}
(N\widetilde{\delta})&=X\left(\lambda_{0}+\frac{\lambda}{2\sqrt{t}}\right)\int_{\lambda_{0}-1/t}^{\lambda_{0}}
\frac{e^{-it\theta(\xi)}f(\xi)}{X_{+}(\xi)\left(\xi-\lambda_{0}-\frac{\lambda}{2\sqrt{t}}\right)}d\xi\notag\\
&+X\left(\lambda_{0}+\frac{\lambda}{2\sqrt{t}}\right)\int_{-\infty}^{\lambda_{0}-1/t}
\frac{e^{-it\theta(\xi)}f_{1}(\xi)}{X_{+}(\xi)\left(\xi-\lambda_{0}-\frac{\lambda}{2\sqrt{t}}\right)}d\xi\notag\\
&+X\left(\lambda_{0}+\frac{\lambda}{2\sqrt{t}}\right)\int_{-\infty}^{\lambda_{0}-1/t}
\frac{e^{-it\theta(\xi)}f_{2}(\xi)}{X_{+}(\xi)\left(\xi-\lambda_{0}-\frac{\lambda}{2\sqrt{t}}\right)}d\xi\notag\\
&=\mathcal {I}_{1}+\mathcal {I}_{2}+\mathcal {I}_{3},
\end{align}
and
\begin{align*}\label{FR-15}
&|\mathcal {I}_{1}|=\lesssim\int_{\lambda_{0}}^{\lambda_{0}-1/t}\frac{f(\xi)}{\left|\xi-\lambda_{0}-\frac{\lambda}{2\sqrt{2t}}\right|}d\xi
\lesssim t^{-l}\left|\log\left|1-\frac{2\sqrt{2}}{\lambda \sqrt{t}}\right|\right|\lesssim t^{-l-1/2},\notag\\
&|\mathcal {I}_{2}|\lesssim\int_{-\infty}^{\lambda_{0}-1/t}\frac{\left|e^{-it\Theta(\xi)}f_{1}(\xi)\right|}{\left|\xi-\lambda_{0}-\lambda/2\sqrt{t}\right|}d\xi
\leq\frac{\sqrt{2}\pi}{2}t^{-l+1}\lesssim t^{-l+1}.
\end{align*}
Following a similar way, as a result of Cauchy's theorem, we can find $\mathcal {I}_{3}$ along the contour $\mathcal {L}_{t}$
take the place of the interval $(-\infty,\lambda_{0}-1/t)$ and get $|\mathcal {I}_{3}|\lesssim t^{-l+1}$.
As a result, it is not hard to see that \eqref{FR-6} holds.

\noindent
\textbf{Note} There exists a similar estimate
\begin{equation}\label{FR-16}
|(N\hat{\delta})(\lambda)|\lesssim t^{-l},~~t\rightarrow\infty,~~\lambda\in\hat{\bar{\mathcal {L}}},
\end{equation}
where
\begin{equation}\label{FR-17}
\hat{\delta}(\lambda)=e^{it\Theta(\lambda)}\left[\delta^{-1}(\lambda)-(\det(\delta))^{-1}(\lambda\mathcal {I})\right]\mathcal {R}^{\dag}(\bar{\lambda}).
\end{equation}

\noindent
\textbf{Theorem 3.11} As $t\rightarrow\infty$
\begin{equation}\label{FR-18}
\textbf{q}(x,t)=\frac{1}{\sqrt{t}}\lim_{\lambda\rightarrow\infty}\left(\lambda M^{0}(x,t;\lambda)\right)_{12}+O\left(\frac{\log t}{\sqrt{t}}\right),
\end{equation}
where $M^{0}(x,t;\lambda)$ meets the RHP
\begin{equation}\label{FR-19}
\left\{ \begin{aligned}
&M_{+}^{0}(\lambda)=M_{-}^{^{0}}(\lambda)J^{0}(\lambda),~~\lambda\in\Sigma^{0},\\
&M^{0}(\lambda)\rightarrow\mathcal {I},~~~~~~~\lambda\rightarrow\infty.
      \end{aligned} \right.
\end{equation}
Here $J^{0}=(\mathcal {I}-\omega^{0}_{+})^{-1}(\mathcal {I}+\omega^{0}_{-})$ and
\begin{equation}\label{FR-10}
\omega^{0}=\omega^{0}_{+}=\left\{ \begin{aligned}
&\left(
  \begin{array}{cc}
    0 & \eta^2\lambda^{2i\nu}e^{-\frac{1}{2}i\lambda^2}\gamma(\lambda_{0}) \\
    0 & 0 \\
  \end{array}
\right),~~\lambda\in\Sigma_{1}^{0},\\
&\left(
  \begin{array}{cc}
    0 & -\eta^2\lambda^{2i\nu}e^{-\frac{1}{2}i\lambda^2}\frac{\gamma(\lambda_{0})}{1+|\gamma(\lambda_{0})|^{2}} \\
    0 & 0 \\
  \end{array}
\right),~~\lambda\in\Sigma_{3}^{0},
      \end{aligned} \right.
\end{equation}
and
\begin{equation}\label{FR-11}
\omega^{0}=\omega^{0}_{+}=\left\{ \begin{aligned}
&\left(
  \begin{array}{cc}
    0 & 0 \\
    \eta^-2\lambda^{-2i\nu}e^{-\frac{1}{2}i\lambda^2}\gamma^{\dag}(\lambda_{0}) & 0 \\
  \end{array}
\right),~~\lambda\in\Sigma_{2}^{0},\\
&\left(
  \begin{array}{cc}
    0 & 0  \\
    -\eta^-2\lambda^{-2i\nu}e^{\frac{1}{2}i\lambda^2}\frac{\gamma^{\dag}(\lambda_{0})}{1+|\gamma(\lambda_{0})|^{2}} & 0 \\
  \end{array}
\right),~~\lambda\in\Sigma_{4}^{0},
      \end{aligned} \right.
\end{equation}
with
\begin{equation}\label{FR-12}
\eta=(4t)^{-\frac{1}{2}i\nu}e^{i\lambda_{0}^{2}t+\widetilde{\chi}(\lambda_{0})}.
\end{equation}

\noindent
\textbf{Proof}
It follows from \eqref{FR-6} and Lemma 3.35 in \cite{PA-1993} that
\begin{equation}\label{FR-13}
\left\|\hat{\omega}-\omega^{0}\right\|_{\mathscr{L}^{\infty}(\Sigma^{0})\bigcap\mathscr{L}^{1}(\Sigma^{0})\bigcap\mathscr{L}^{2}(\Sigma^{0})}\lesssim
\frac{\log t}{\sqrt{t}}.
\end{equation}
As a result,
\begin{align}\label{FR-14}
\int_{\Sigma'}\left(\left(1-C_{\omega'}\right)^{-1}\mathcal {I}\right)(\xi)\omega'(\xi)d\xi
&=\int_{\Sigma'}\left(N^{-1}\left(1-C_{\omega'}\right)^{-1}N\mathcal {I}\right)(\xi)\omega'(\xi)d\xi\notag\\
&=\int_{\Sigma'}((1-C_{\hat{\omega}})^{-1}\mathcal {I})\left(2(\xi-\lambda_{0})\sqrt{t}\right)N\omega'\left(2(\xi-\lambda_{0})\sqrt{t}\right)d\xi\notag\\
&=\frac{1}{2\sqrt{t}}\int_{\Sigma_{0}}\left(\left(1-C_{\hat{\omega}}\right)^{-1}\mathcal {I}\right)(\xi)\hat{\omega}(\xi)d\xi \notag\\
&=\frac{1}{2\sqrt{t}}\int_{\Sigma_{0}}\left(\left(1-C_{\omega^{0}}\right)^{-1}\mathcal {I}\right)(\xi)\omega^{0}(\xi)d\xi+O\left(\frac{\log t}{t}\right).
\end{align}
For $\lambda\in\mathbb{C}\setminus\Sigma^{0}$, let
\begin{equation}\label{FR-15}
M^{0}(\lambda)=\mathcal {I}+\frac{1}{2\pi i}\int_{\Sigma_{0}}\frac{\left(\left(1-C_{\omega^{0}}\right)^{-1}\mathcal {I}\right)(\xi)\omega^{0}(\xi)}{\xi-\lambda}d\xi.
\end{equation}
Then $M^{0}$ can deal with the above RHP.
From above expressions and Lemma 3.8, it is straightforward to derive this theorem.

\noindent
\textbf{Note} Particularly, if
\begin{equation}\label{FR-16}
M^{0}(\lambda)=\mathcal {I}+\frac{M_{1}^{0}}{\lambda}+O\left(\frac{1}{\lambda^{2}}\right),~~\lambda\rightarrow\infty,
\end{equation}
then
\begin{equation}\label{FR-17}
\textbf{q}(x,t)=\frac{1}{\sqrt{t}}\left(M_{1}^{0}\right)_{12}+O\left(\frac{\log t}{t}\right).
\end{equation}

\subsection{Solving the Model Problem}
In order to give $M_{1}^{0}$ explicitly, it is worth considering the following transformation
\begin{equation}\label{SMP-1}
\Psi(\lambda)=H(\lambda)\lambda^{-iv\sigma}e^{i\lambda^2\sigma/4},~~H(\lambda)=\eta^{\sigma}M^{0}(\lambda)\eta^{-\sigma},
\end{equation}
which indicates that
\begin{equation}\label{SMP-2}
\Psi_{+}(\lambda)=\Psi_{-}(\lambda)v(\lambda_{0}),~~v=\lambda^{iv\hat{\sigma}}e^{-i\lambda^2\hat{\sigma}/4}\eta^{\hat{\sigma}}J^{0}.
\end{equation}
As the jump matrix is constant along each ray, we have
\begin{equation}\label{SMP-3}
\frac{d\Psi_{+}(\lambda)}{d\lambda}=\frac{d\Psi_{-}(\lambda)}{d\lambda}v(\lambda_{0}),
\end{equation}
from which it follows that $\frac{d\Psi_{-}(\lambda)}{d\lambda}\Psi^{-1}(\lambda)$
has no jump discontinuity along any of the rays.
Additionally, from the relation between $H(\lambda)$ and $\Psi(\lambda)$, we have
\begin{align}\label{SMP-4}
\frac{d\Psi(\lambda)}{d\lambda}\Psi^{-1}(\lambda)&=\frac{dH(\lambda)}{d\lambda}H^{-1}(\lambda)+\frac{1}{2}i\lambda H(\lambda)\sigma H^{-1}(\lambda)
-\frac{i\nu}{\lambda}H(\lambda)\sigma H^{-1}(\lambda)\notag\\
&=O\left(\frac{1}{\lambda}\right)+\frac{1}{2}i\lambda\sigma-\frac{1}{2}i\eta^{\sigma}\left[\sigma,M_{1}^{0}\right]\eta^{-\sigma}.
\end{align}
It follows from the Liouville's theorem that
\begin{equation}\label{SMP-5}
\frac{d\Psi(\lambda)}{d\lambda}=\frac{1}{2}i\lambda\sigma\Psi(\lambda)+\beta\Psi(\lambda),
\end{equation}
where
\begin{equation}\label{SMP-6}
\beta=-\frac{1}{2}i\eta^{\sigma}\left[\sigma,M_{1}^{0}\right]\eta^{-\sigma}=\left(
                                                                   \begin{array}{cc}
                                                                     0 & \beta_{12} \\
                                                                     \beta_{21} & 0 \\
                                                                   \end{array}
                                                                 \right).
\end{equation}
Particularly,
\begin{equation}\label{SMP-7}
\left(M_{1}^{0}\right)=-i\eta^{2}\beta_{12}.
\end{equation}
It is further possible to find that the solution of the RHP for $M^{0}(\lambda)$ is unique,
and therefore we have an identity
\begin{equation}\label{SMP-8}
\left(M^{0}(\bar{\lambda})\right)^{\dag}=\left(M^{0}(\lambda)\right)^{-1},
\end{equation}
which implies that $\beta_{12}=-\beta_{21}^{\dag}$.
Let
\begin{equation*}\label{SMP-9}
\Psi(\lambda)=\left(
                \begin{array}{cc}
                  \Psi_{11}(\lambda) & \Psi_{12}(\lambda) \\
                  \Psi_{21}(\lambda) & \Psi_{22}(\lambda) \\
                \end{array}
              \right).
\end{equation*}
From \eqref{SMP-5} we obtain
\begin{align}\label{SMP-10}
&\frac{d^{2}\Psi_{11}(\lambda)}{d\lambda^{2}}=\left(\beta_{12}\beta_{21}-\frac{i}{2}-\frac{\lambda^2}{4}\right)\Psi_{11},\notag\\
&\beta_{21}\Psi_{21}(\lambda)=\frac{d\Psi_{11}(\lambda)}{d\lambda}+\frac{i}{2}\lambda\Psi_{11},\notag\\
&\frac{d^{2}\beta_{12}\Psi_{22}(\lambda)}{d\lambda^2}=\left(\beta_{12}\beta_{21}+\frac{i}{2}-\frac{\lambda^2}{4}\right)\beta_{12}\Psi_{22}(\lambda),\notag\\
&\Psi_{12}(\lambda)=\frac{1}{\beta_{12}\beta_{21}}\left(\frac{d\beta_{12}\Psi_{22}(\lambda)}{d\lambda}-\frac{i}{2}\lambda\beta_{12}\Psi_{22}(\lambda)\right).
\end{align}
As we all know, the Weber's equation
\begin{equation*}\label{SMP-11}
\frac{d^{2}g(\zeta)}{d\zeta^2}+\left(a+\frac{1}{2}-\frac{\zeta^2}{4}\right)g(\zeta)=0,
\end{equation*}
admits the solution
\begin{equation*}\label{SMP-12}
g(\zeta)=c_{1}D_{a}(\zeta)+c_{2}D_{a}(-\zeta),
\end{equation*}
where $D_{a}(\cdot)$ represents the standard parabolic-cylinder function and satisfies
\begin{align}\label{SMP-13}
&\frac{dD_{a}(\zeta)}{d\zeta}+\frac{\zeta}{2}D_{a}(\zeta)-aD_{a-1}(\zeta)=0,\notag\\
&D_{a}(\pm\zeta)=\frac{\Gamma(1+a)e^{i\pi a/2}}{\sqrt{2\pi}}D_{-a-1}(\pm i\zeta)+\frac{\Gamma(1+a)e^{-i\pi a/2}}{\sqrt{2\pi}}D_{-a-1}(\mp i\zeta).
\end{align}
As $\zeta\rightarrow\infty$, from \cite{MJ-2003} we have
\begin{align}\label{SMP-14}
&D_{a}(\zeta)=\notag\\
&\left\{ \begin{aligned}
&\zeta^{a}e^{-\zeta^2/4}\left(1+O\left(\frac{1}{\zeta^{2}}\right)\right),~~|\mbox{arg}\zeta|<\frac{3\pi}{4},\\
&\zeta^{a}e^{-\zeta^2/4}\left(1+O\left(\frac{1}{\zeta^{2}}\right)\right)-\frac{\sqrt{2\pi}}{\Gamma(-a)}e^{a\pi i+\zeta^2/4}\left(1+O\left(\frac{1}{\zeta^{2}}\right)\right),~~\frac{3\pi}{4}<\mbox{arg}\zeta<\frac{5\pi}{4},\\
&\zeta^{a}e^{-\zeta^2/4}\left(1+O\left(\frac{1}{\zeta^{2}}\right)\right)-\frac{\sqrt{2\pi}}{\Gamma(-a)}e^{-a\pi i+\zeta^2/4}\left(1+O\left(\frac{1}{\zeta^{2}}\right)\right),~~-\frac{5\pi}{4}<\mbox{arg}\zeta<-\frac{3\pi}{4},
      \end{aligned} \right.
\end{align}
where $\Gamma$ is the Gamma function.

Choose $a=i\beta_{12}\beta_{21}$, we have
\begin{align*}\label{SMP-15}
&\Psi_{11}(\lambda)=c_{1} D_{a}\left(e^{-3\pi i/4}\lambda\right)+c_{2}D_{a}\left(e^{\pi i/4}\lambda\right),\notag\\
&\beta_{12}\Psi_{22}(\lambda)=c_{3}D_{-a}\left(e^{3\pi i/4}\lambda\right)+c_{2}D_{-a}\left(e^{-\pi i/4}\lambda\right),
\end{align*}
As $\mbox{arg}\lambda\in(-\frac{\pi}{4},\frac{\pi}{4})$ and $\lambda\rightarrow\infty$, we arrive at
\begin{equation*}\label{SMP-16}
\Psi_{11}(\lambda)\lambda^{-i\nu}e^{i\lambda^2/4}\rightarrow \mathcal {I},~~\Psi_{22}(\lambda)\lambda^{i\nu}e^{-i\lambda^2/4}\rightarrow 1.
\end{equation*}
Then
\begin{equation*}\label{SMP-17}
\Psi_{11}(\lambda)=e^{\pi\nu/4}D_{a}\left(e^{\pi i/4}\lambda\right),~~\nu=\beta_{21}\beta_{12},
~~\beta_{12}\Psi_{22}(\lambda)=\beta_{12}e^{\pi\nu/4}D_{-a}\left(e^{-\pi i/4}\lambda\right).
\end{equation*}
Consequently,
\begin{equation*}\label{SMP-18}
\beta_{12}\Psi_{21}(\lambda)=ae^{\pi(\nu+i)/4}D_{a-1}\left(e^{\pi i/4}\lambda\right),~~
\Psi_{21}(\lambda)=\beta_{12}e^{\pi(\nu-3i)/4}D_{-a-1}\left(e^{-\pi i/4}\lambda\right).
\end{equation*}
For $\mbox{arg}\lambda\in(-\frac{3\pi}{4},-\frac{\pi}{4})$ and $\lambda\rightarrow\infty$, we have
\begin{equation*}\label{SMP-19}
\Psi_{11}(\lambda)\lambda^{-i\nu}e^{i\lambda^2/4}\rightarrow\mathcal {I},~~\Psi_{22}(\lambda)\lambda^{i\nu}e^{-i\lambda^2/4}\rightarrow1.
\end{equation*}
Then
\begin{equation*}\label{SMP-20}
\Psi_{11}(\lambda)=e^{\pi\nu/4}D_{a}\left(e^{\pi i/4}\lambda\right),
~~\beta_{12}\Psi_{22}(\lambda)=\beta_{12}e^{-3\pi\nu/4}D_{-a}\left(e^{3\pi i/4}\lambda\right).
\end{equation*}
Hence,
\begin{equation*}\label{SMP-21}
\beta_{21}\Psi_{21}(\lambda)=ae^{\pi(\nu+i)/4}D_{a-1}\left(e^{\pi i/4}\lambda\right),~~
\Psi_{12}(\lambda)=\beta_{12}e^{\pi(\nu-3i)/4}D_{-a-1}\left(e^{-3\pi i/4}\lambda\right).
\end{equation*}
Following the ray $\mbox{arg}\lambda=-\frac{\pi}{4}$,
\begin{align*}\label{SMP-22}
&\Psi_{+}(\lambda)=\Psi_{-}(\lambda)\left(
                                     \begin{array}{cc}
                                        1 & \gamma(\lambda_{0}) \\
                                        0 & \mathcal {I} \\
                                     \end{array}
                                   \right),\notag\\
&\beta_{12}e^{\pi(i-3\nu)/4}D_{-a-1}\left(e^{3\pi i/4}\lambda\right)=e^{\pi\nu/4}D_{a}\left(e^{\pi i/4}\lambda\right)\gamma(\lambda_{0})
+\beta_{12}e^{\pi(v-3i)/4}D_{-a-1}\left(e^{-\frac{\pi i}{4}}\lambda\right).
\end{align*}
It follows from \eqref{SMP-13} that
\begin{equation*}\label{SMP-23}
D_{a}\left(e^{\pi i/4}\lambda\right)=\Gamma(1+a)e^{i\pi a/2}D_{-a-1}\left(e^{3\pi i/4}\lambda\right)+
\frac{\Gamma(1+a)e^{-i\pi a/2}}{\sqrt{2\pi}}D_{-a-1}\left(e^{-\pi i/4}\lambda\right).
\end{equation*}
We then analyze the coefficients of the two independent functions and get
\begin{equation}\label{SMP-25}
\beta_{12}=\frac{e^{\pi\nu/2-\pi i/4}\Gamma(1+a)}{\sqrt{2\pi}}\gamma(\lambda_{0})=
\frac{e^{\pi\nu/2+\pi i/4}\nu\Gamma(i\nu)}{\sqrt{2\pi}}\gamma(\lambda_{0}).
\end{equation}
Summarizing the above results, the following Theorem 3.12 can be easily established.

\noindent
\textbf{Theorem 3.12}
Suppose that $(q_{1},q_{2},q_{3})$ possesses the solution for the Cauchy problem of the 3-component Manakov system \eqref{tc-NLS}
with $q_{1,0},q_{2,0},q_{3,0}\in\mathcal {G}$. Then $|\frac{x}{t}|\leq C$,
the leading asymptotics of $(q_{1},q_{2},q_{3})$ admits the following explicit form
\begin{equation}\label{SMP-26}
\textbf{q}(x,t)=(q_{1},q_{2},q_{3})=\frac{\nu}{2\sqrt{\pi t}}\Gamma(i\nu)(4t)^{i\nu}\gamma(\lambda_{0})e^{2i\lambda_{0}^2t+2\widetilde{\chi}(\lambda_{0})+\frac{\pi\nu}{2}-\frac{\pi i}{4}}
+O\left(\frac{\log t}{t}\right),
\end{equation}
where $\lambda_{0}=-\frac{x}{2t}$, $C$ is a constant, $\Gamma(\cdot)$
represents a Gamma function, the vector-valued function $\gamma(\lambda)$ is expressed by \eqref{TH-3}.
Additionally, $\nu$ and $\widetilde{\chi}(\lambda_{0})$ are determined by \eqref{TRH-5}.

$~~~~~~~~~~~~~~~~~~~$
{\rotatebox{0}{\includegraphics[width=8.0cm,height=6.5cm,angle=0]{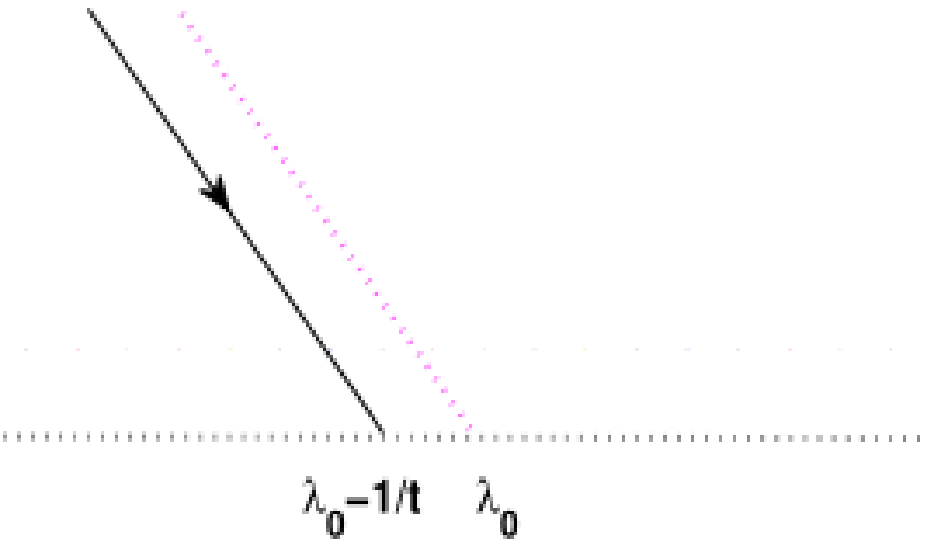}}}\\
\noindent {\small \textbf{Figure 4.} The contour $\mathcal {L}_{t}$ is defined by \eqref{FR-14}. \\}

\section{Conclusions}

In this paper, by generalizing the nonlinear steepest descent method,
we have studied the long-time behavior for the Cauchy problem of the 3-component Manakov system  \eqref{tc-NLS} with a $4\times4$ Lax pair.
The nonlinear steepest descent technique
is good for exploring the long-time asymptotics of integrable equations with Lax pairs and is also good for
integrable systems with higher-order spectral problems.
Here there are two reasons for choosing the 3-component Manakov system  \eqref{tc-NLS} as a model problem:
(I) due to its physical interest.
To well model important types of nonlinear physical
phenomena in a proper way, there is a necessity to go
beyond the standard NLS description.
A crucial development consists of the investigation of coupled nonlinear models, as many physical systems comprise
interacting wave components of distinct modes, polarizations or frequencies. In recent years, the multi-components coupled NLS
equations have become a topic of intense research in the
field of mathematical physics, since the components are
usually more than one practically for many physical phenomena.
(II) There have been the asymptotics for several nonlinear integrable equations with $2\times2$ and $3\times3$
Lax pairs, for instance, KdV equation, NLS equation, mKdV equation, Camassa-Holm equation, derivative NLS equation,
Sasa-Satsuma equation, coupled NLS equations etc.
However, but there is just a little of literature about nonlinear integrable equations with $4\times4$ Lax pairs.
As a result, how to discuss the asymptotic behavior of nonlinear integrable equations with $4\times4$ Lax pairs is interesting and meaningful.
In \cite{PA-1993},
the function $\delta$ can be solved explicitly by the Plemelj formula because $\delta$ meets a scalar RHP.
However, the function $\delta$ in our work admits $3\times3$ matrix RHP.
The unsolvability of the $3\times3$ matrix function $\delta$ is a challenge for us.
Noticing that our
purpose is to investigate the asymptotic behavior for solution of the 3-component Manakov system \eqref{tc-NLS},
a natural idea is using the available function $\det \delta$ to approximate $\delta$ by error control.

Finally, we remark that the nonlinear steepest descent technique
can be used to analyze the Long-time asymptotics of many nonlinear integrable equations with nonzero boundary conditions.
Thus,
it is essential  to discuss whether
the Long-time asymptotics of the 3-component Manakov system with nonzero boundary conditions can be obtained by using
the nonlinear steepest descent technique? These will be left for future discussions.

\section*{Acknowledgements}
\hspace{0.3cm}
We express our sincere thanks to
the editor and reviewers for their valuable comments.
This work is supported by the National Key Research and Development Program of
China under Grant No. 2017YFB0202901 and the National Natural Science Foundation of China under Grant No.11871180.

\end{document}